\documentclass[11pt]{article}
\usepackage{amsmath,latexsym,amssymb,graphicx,epsf,amsfonts,amsthm}
\usepackage[margin=4.46cm]{geometry}
\usepackage{epsf}
\usepackage{color}
\usepackage{mathrsfs}
\usepackage{epstopdf}

\numberwithin{equation}{section} 

\newtheorem{thm}{Theorem}
\newtheorem{theorem}{Theorem}

\newtheorem{alg}[thm]{Algorithm}

\newtheorem{ass}{Assumption}

\newtheorem{rem}[thm]{Remark}

\begin{document}
\title{A Geometric Blind Source Separation \\Method Based on Facet Component Analysis}
\author{Penghang Yin\thanks{Department of Mathematics, University of California at Irvine, Irvine, CA 92697, USA.}, Yuanchang Sun\thanks{Department of Mathematics and Statistics, Florida International University, Miami, FL 33199, USA.}, and Jack Xin$^{*}$}
\date{}
\maketitle

\begin{abstract}
Given a set of mixtures, blind source separation attempts to retrieve the source signals without or with very little information of the the mixing process.  We present a geometric approach for blind separation of nonnegative linear mixtures termed {\em facet component analysis} (FCA).  The approach is based on facet identification of the underlying cone structure of the data.
Earlier works focus on recovering the
cone by locating its vertices (vertex component analysis or VCA) based on a mutual sparsity
condition which requires each source signal to possess a stand-alone peak in its spectrum.
We formulate alternative conditions so that enough data points fall on the facets of a cone
instead of accumulating around the vertices.
To find a regime of unique solvability, we make use of both geometric and
density properties of the data points, and develop an efficient facet identification method by combining data classification and linear regression.
For noisy data, we show that denoising methods may be employed, such as
the total variation technique in imaging processing, and principle component analysis.  We show computational results on
nuclear magnetic resonance spectroscopic data to substantiate our method.

\end{abstract}
\thispagestyle{empty}
\newpage
\setcounter{equation}{0} \setcounter{page}{1}
\section{Introduction}

Blind source separation (BSS) is a major area of research in signal and image processing \cite{Cic}.
It aims at recovering source signals from their mixtures with minimal knowledge of the mixing environment.
The applications of BSS range from engineering to neuroscience.
A recent emerging research direction of BSS is to identify chemical explosives
and biological agents from their spectral sensing mixtures
recorded by various spectroscopy such as Nuclear Magnetic Resonance (NMR),
Raman spectroscopy, Ion-mobility spectroscopy (IMS), and
differential optical absorption spectroscopy (DOAS),etc.
Spectral sensing is a critical area in national security and a vibrant scientific area.
The advances of modern imaging and spectroscopic technology have made it possible to
classify pure chemicals by their spectral features.
However, mixtures of chemicals subject to changing background and environmental noise pose
additional challenges.  The goal of this paper is to develop
a BSS method to process data in the presence of noise based on geometric spectral properties.

To separate the spectral mixtures, one needs to solve the
following matrix decomposition problem
\begin{equation}
\label{bss}
X = A\,S + N\;,
\end{equation}
where $A \in \mathbb{R}^{m\times n}$ is a full rank unknown basis (dictionary) matrix
or the so called mixing matrix in some applications, $N \in \mathbb{R}^{m\times p}$
is an unknown noise matrix,
$S = [s(1),\dots, s(p)] \in  \mathbb{R}^{n\times p}$ is the unknown source matrix containing
signal spectra in its rows. Here $p$ is the number of data samples,
$m$ is the number of observations, and $n$ is the number of sources.
Various BSS methods have been proposed based on {\em a priori} knowledge of source signals
such as statistical independence, sparseness, nonnegativity,
\cite{Bof,Choi,Cic,Hoyer,Hyv,LX_08,LX_09a,LX_09b,NN05,SX1,SX2} among others.
As a matrix factorization problem in the noise free case ($N=0$),
BSS has permutation and scaling ambiguities in its
solutions similar to factorizing a large number into product of primes.
For any permutation matrix $P$ and invertible diagonal matrix $\Lambda$,
($AP\Lambda$, $\Lambda^{-1}P^{-1}S$) is another pair equivalent to the solution $(A,S)$,
since
\begin{equation}
\label{eqvbss} X = A\,S = (A\,P\Lambda)(\Lambda^{-1}P^{-1}S).
\end{equation}

\begin{figure}
\includegraphics[height=7cm,width=7cm]{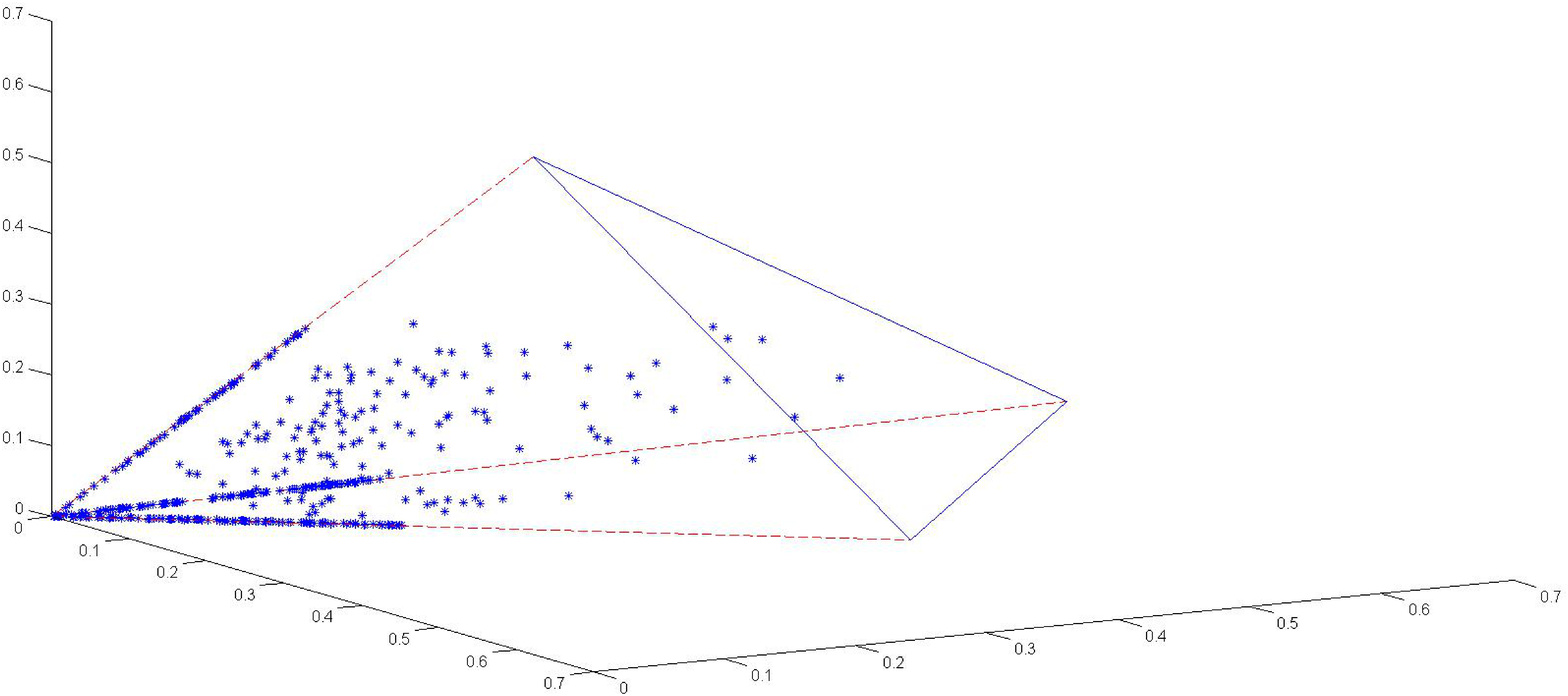}
\includegraphics[height=7cm,width=7cm]{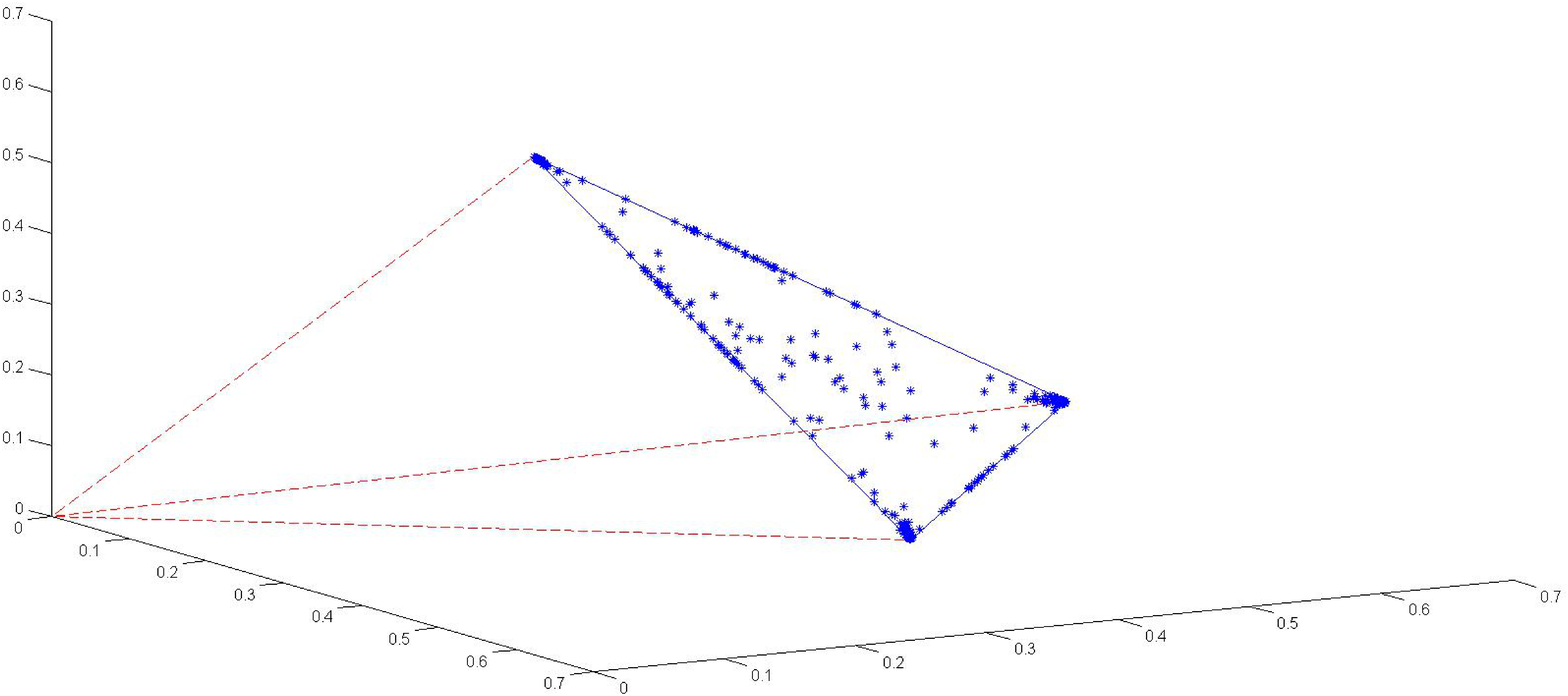}
\caption{A cloud of data points (columns of $X$), left plot, is rescaled to lie
on a plane determined by three vertices of the cone (right plot). }
 \label{cone1}
\end{figure}

Recently there has been active research on BSS by exploiting data geometry
\cite{Boardman_93,Bobin_1,Chang_07,Kli,NN05,NN10,Nas,SX1, SX2,SUN_XIN_JSC, SUN_XIN_siim, Winter_99}.
For simplicity, let $N=0$ in (\ref{bss}). The geometric observation \cite{Boardman_93,NN05,Winter_99}
is that if each row of $S$ has a dominant peak at some location (column number)
where other rows have zero elements, then the problem of finding the columns of the
mixing matrix $A$ reduces to the identification of the edges of a minimal cone
containing the columns of mixture matrix $X$. In hyperspectral imaging (HSI), the
condition is known as pixel purity assumption (PPA,\cite{Chang_07}). In other words, each pure material of interest exists by itself somewhere on the ground.
The PPA based convex cone method (known as N-findr \cite{Winter_99}) is now a benchmark in HSI, see \cite{Chang_07,NN05,NN10, Nas, SX1} for its more recent variants.  The method termed {\em vertex component analysis} (VGA) proposed in \cite{Nas} is worth mentioning here being a fast unmixing algorithm for hyperspectral data.
In Nuclear Magnetic Resonance (NMR) spectroscopy which motivates our work here,
PPA was reformulated by Naanaa and Nuzillard in \cite{NN05}.
The source signals are only required to be non-overlapping at some locations of
acquisition variable (e.g. frequency) which leads
to a dramatic mathematical simplification of a general non-negative
matrix factorization problem (\ref{bss}) which is non-convex \cite{Lee}.
More precisely, the source matrix $S\geq 0$ is assumed to satisfy the following
\begin{ass}[NNA]: For each $i\in\{1,2,\dots,n \}$ there exists an $j_{i}\in
\{1,2,\dots,p\}$ such that $s_{i,j_i}>0$ and $s_{k,j_i}=0\; (k=1,\dots,i-1,i+1,\dots,n)\;.$
\end{ass}
Simply put, the stand-alone peaks possessed by each source allow formation
of a convex cone enclosing all the columns of $X$, and
the edges (or vertices) of the cone are the columns of the mixing matrix.
To illustrate the idea, let us consider an example of three sources and three
mixtures ($A \in \mathbb{R}^{3\times 3}, X \in \mathbb{R}^{3\times p}, S\in \mathbb{R}^{3\times p}$,
and $A$ is non-singular).
\begin{eqnarray*}
& & \bigl [X(:,1),X(:,2),\cdots,X(:,p)\bigr] \\
 & & = \bigl[A(:,1),A(:,2),A(:,3)\bigr]\cdot
\left(
\begin{array}{cccccccccc}
        * & \cdots & * & {\bf 1 }     & {0}      & {0}     & * & \cdots & *\\
        * & \cdots & * & {0 }         & {\bf 1}  & {0}     & * & \cdots & *\\
        * & \cdots & * & {0 }         & {0}      & {\bf 1} & * & \cdots & *
       \end{array}
 \right)\;,
\end{eqnarray*}
where $X(k,:)$ represents the $k$-th row of $X$, ${\bf 1}$ represents nonzero entry, and the three ${\bf 1}'s$ are three stand-alone peaks. It can be seen that $A$'s columns are actually among those of $X$'s up to constants, and other columns of $X$ are
nonnegative linear combinations of columns of $A$.
Geometrically, the columns of $A$ span a convex cone enclosing all $X$'s columns.
The estimation of $A$ is equivalent to the identification of this cone.
As a matter of fact, the vertices of the cone are the columns of $A$.
Fig. \ref{cone1} shows this vertex-represented cone containing all data points.
To find the vertices (or edges) of the cone, the following optimization problem
is solved for each column $k$
of $X$:
\begin{eqnarray*}
\label{LPNF}
c &  = & \min \sum^p_{j = 1, j\neq k} \lambda_j \\
&\mathrm{s.t.}&  \;\; \sum^{p}_{j = 1, j\neq k}X(:,j) \lambda_j =
X(:,k)\;\;, \lambda_j\geq 0.
\end{eqnarray*}
It is shown \cite{Dul} that {\it $X(:,k)$ is an edge of the convex
cone if and only if the optimal objective function value $c^*$ is greater than 1.}

If the data are contaminated by noise, the following optimization problem
is solved for each $k$ \cite{NN05}:
{\allowdisplaybreaks
\begin{equation}
\label{LPNP} \min\; \mathrm{score} = \min_{\lambda_j \geq 0} \; \|\sum^{p}_{j = 1, j\neq k}X(:,j) \lambda_j - X(:,k)
 \|^2_2.
\end{equation}}
A score is associated with each column of $X$.  Columns with low scores are
unlikely to be a column of $A$ because this
column is roughly a nonnegative linear combination of the other columns of $X$.
On the other hand, a high score means that the
corresponding column is far from being a nonnegative linear combination of other columns.
 The $n$ rows from $X$ with highest scores are selected to form $A$, the mixing matrix.

\begin{figure}
\includegraphics[height=7cm,width=7cm]{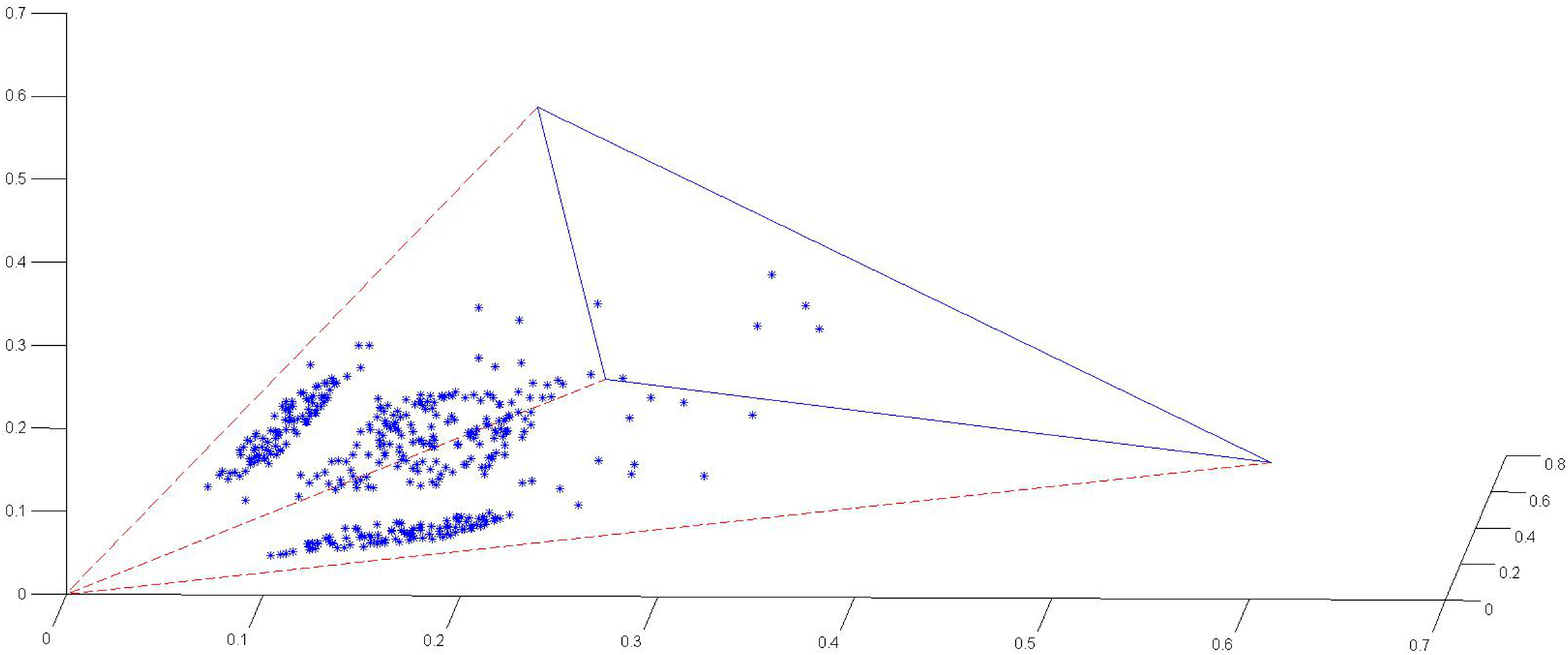}
\includegraphics[height=7cm,width=7cm]{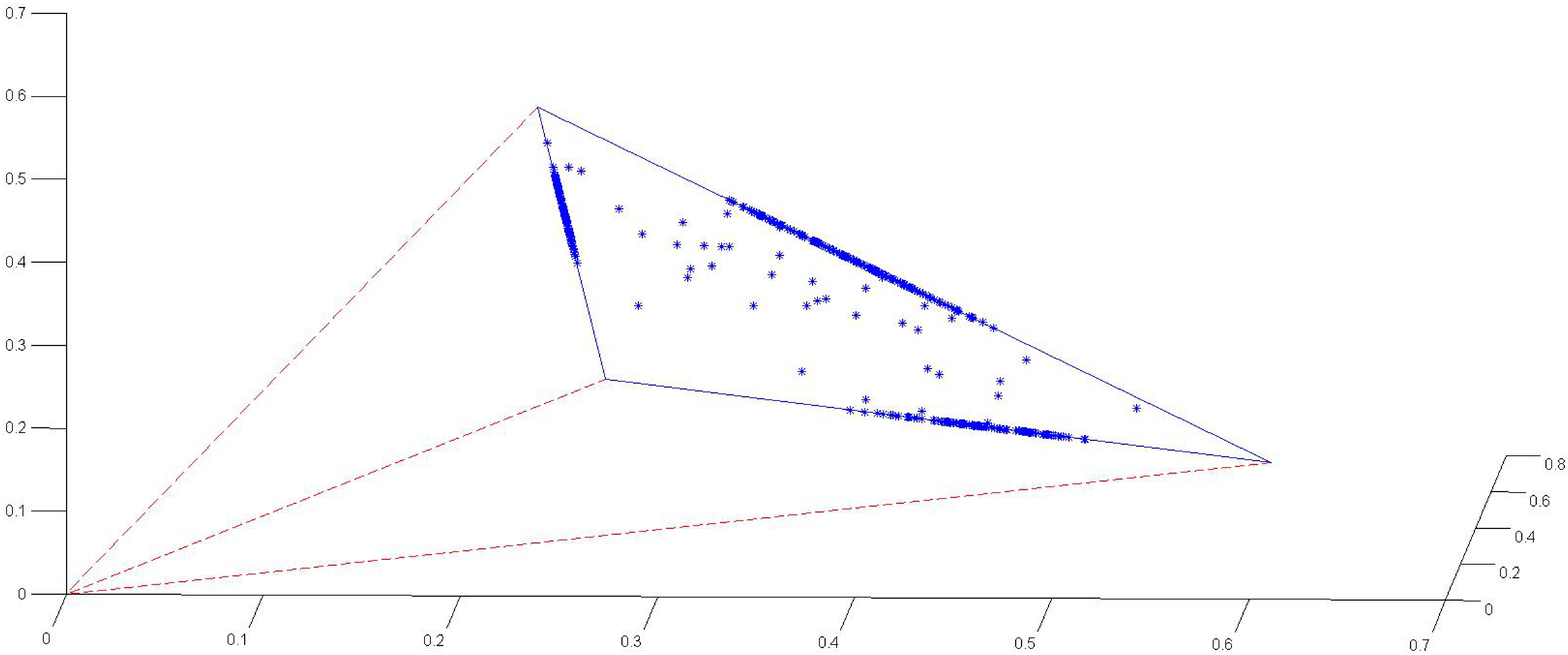}
\caption{Scatter plots of columns of $X$. The cloud of data point (left) rescaled to lie
on a plane determined by the three vertices of the cone (right). }
 \label{cone_facet}
\end{figure}

Though the vertex based convex cone method above is geometrically elegant,
the working condition is still restrictive.
The success of the cone method highly depends on the recognition of its vertices.
If PPA or NNA is violated, the vertices are not in the data matrix $X$, and
may not be the primary objects for identification.  In this paper, we consider such a scenario
where data points (scaled columns of $X$) lie either inside or on the facets of a convex cone,
yet none of them are located on edges (vertices), see Fig. \ref{cone_facet} for an example.
The cone structure will be reconstructed from its facets instead of the vertices.  A facet component analysis (FCA) should be pursued for the data considered.
The appropriate source condition for this case is
that {\em most columns of source matrix $S$ ($m$ rows) possess $m-1$ nonzero entries}.
We shall pursue a more precise study on the source condition for solvability in the next section.
The problem above calls for an identification method of flat submanifolds from a point cloud,
as the facets of a convex cone in general lie on hyperplanes.
 The vertices are obtained from the intersections of the hyperplanes expanded from the facets.
We shall study the identification of these flat manifolds and the subequent
recovery of the source signals. The extraction of
meaningful geometric data structures help the data matrix factorization
and reduce the computational cost.

Recently, a dual cone based approach to BSS problem was proposed in \cite{NN10}.
The method first calculates the dual of the data cone, then selects a subset of the
source signals from a set of source candidates by means of a greedy process.
Specifically, the first step consists in computing a generating matrix of
the dual of the data cone by the double description method \cite{DDM}.
The second step consists in extracting an
estimate of the source matrix from a larger matrix which is the product of the
generating matrix with its transpose.  Multiple solutions can be obtained in the second step,
so the author imposed an additional orthogonality constraint on the source sigals to get
 a sub-optimal solution.  Overall, the method proposed in \cite{NN10} appears indirect
and computationally unwieldy, besides requiring the orthogonality of source signals.
Here we opt to solve the problem directly by identifying the
facets of the data cone under a unique solvability condition.

The rest of the paper is organized as follows.  In section 2, we propose a new source condition for
solving (\ref{eqvbss}) based on the geometric structure of the data points,
and develop an algorithm for facet
identification and reconstruction of the associated convex cone.
In section 3, we present computational examples and results. For noisy data, we show that a denoising method maybe needed to help remove or reduce the noise.  A denoising method based on total variation and distance function is discussed.
Concluding remarks and future work are in Section 4.

The work was partially supported by NSF-ATD grant DMS-0911277.

\section{Proposed Method}

Let us consider the noiseless case in order to illustrate the basic idea behind our method.
The first source condition on the problem is as follows:
\begin{ass}\label{nonNN}
$S$ contains at least $m-1$ linearly independent column vectors orthogonal to each unit coordinate vector
$e_{i}$, $i = 1,\ldots,m$.
\end{ass}
For the case $m = 3$, we have three mixtures and three sources (i.e, $A \in \mathbb{R}^{3\times 3}, X \in \mathbb{R}^{3\times p}, S\in \mathbb{R}^{3\times p}$, and $A$ is non-singular), one example can be:
\begin{eqnarray*}
& & \bigl [X(:,1),X(:,2),\cdots,X(:,p)\bigr]\\
 & & = \bigl[A(:,1),A(:,2),A(:,3)\bigr]\cdot
\left(
\begin{array}{cccccccccccc}
        * & \cdots & * & {\bf 1 } & {\bf 2 }  & {0}      & {0}      & {\bf 1} & {\bf 2} & * & \cdots & *\\
        * & \cdots & * & {\bf 2 } & {\bf 1 }  & {\bf 1}  & {\bf 2}  & {0}     & {0}     & * & \cdots & *\\
        * & \cdots & * & {0}      & {0}       & {\bf 2}  & {\bf 1}  & {\bf 2} & {\bf 1} & * & \cdots & *
       \end{array}
 \right)\;,
\end{eqnarray*}
We note in passing that PPA or NNA (stand-alone peak) assumption
actually is a special (more restrictive) case of the assumption \ref{nonNN} above.

Let $M$ be a matrix of $\mathbb{R}^{m\times p}$ , the subset of $\mathbb{R}^{m}$ defined by
$$
\mathcal{M} = cone(M) = \{ M\,\alpha\,| \alpha \succcurlyeq 0 \}
$$
is a convex cone, and $M$ is said to be a
generating matrix of $\mathcal{M}$ since every element of $\mathcal{M}$ is a nonnegative
linear combination of $M$ columns.  Let $\mathcal{X}$ = $cone(X)$ and $\mathcal{A}$ = $cone(A)$.
Under Assumption \ref{nonNN}, we have the following (or combining Lemma 3 and Lemma 5 of \cite{NN10}):

\begin{theorem}\label{thm 1}
If $X = A\,S$, and $A$, $S \succcurlyeq 0$, then $\mathcal{X} \subseteq \mathcal{A}$. Moreover, each facet of $\mathcal{A}$ contains a facet of $\mathcal{X}$.
\end{theorem}
For readers' convenience, a short proof is given below.

\begin{proof}
$\forall \, x \in \mathcal{X}$, let $x = X\,\alpha, \alpha \succcurlyeq 0$.
So $x = A\,S\,\alpha = A\,(S\,\alpha)$, where $S\,\alpha \succcurlyeq 0$. So clearly $x \in \mathcal{A}$.

The second claim follows as we notice that: (1) $\mathcal{A}$ has $m$ facets and
each one is spanned by $m-1$ column vectors of $A$; (2) Using Assumption \ref{nonNN},
$X = A\,S$ has at least $m-1$ linearly independent column vectors located in each facet of $\mathcal{A}$;
(3) $\mathcal{X} \subseteq \mathcal{A}$.
\end{proof}

Let $\textbf{x} = (x_{1},\dots,x_{m})^{T}$, $\textbf{1} = (1,\dots,1)^{T}$.
Since $A$ is nonsingular, $\mathcal{A}$ has $m$ extreme directions and thus has $\binom{m}{m-1} = m$ facets.
Based on Theorem \ref{thm 1}, our method aims to identify the $m$ facets of
$\mathcal{X}$ contained in the facets of $\mathcal{A}$.
If we project $X$'s column vectors onto the plane $\textbf{x}^{T}\cdot\textbf{1} = 1$ (i.e. $L_{1}$-normalize
the column vectors), the resulting data points together with the
origin $\textbf{O} = (0,\ldots,0)^{'}$ has a $m$-dimensional convex hull
with $l$ facets($l \geq m+1$). Let us denote it by $Conv(X)$. In fact, $Conv(X)$ results from the
cone $\mathcal{X}$ being truncated by the plane $\textbf{x}^{T}\cdot\textbf{1} = 1$.
We then acquire all the facets and the associated vertices of $Conv(X)$.
This can be done by means of the function 'Convhulln' from the Matlab library.
It is based on 'Qhull' which implements the Quickhull algorithm for computing the convex hull.

\begin{figure}
\includegraphics[height=7cm,width=15cm]{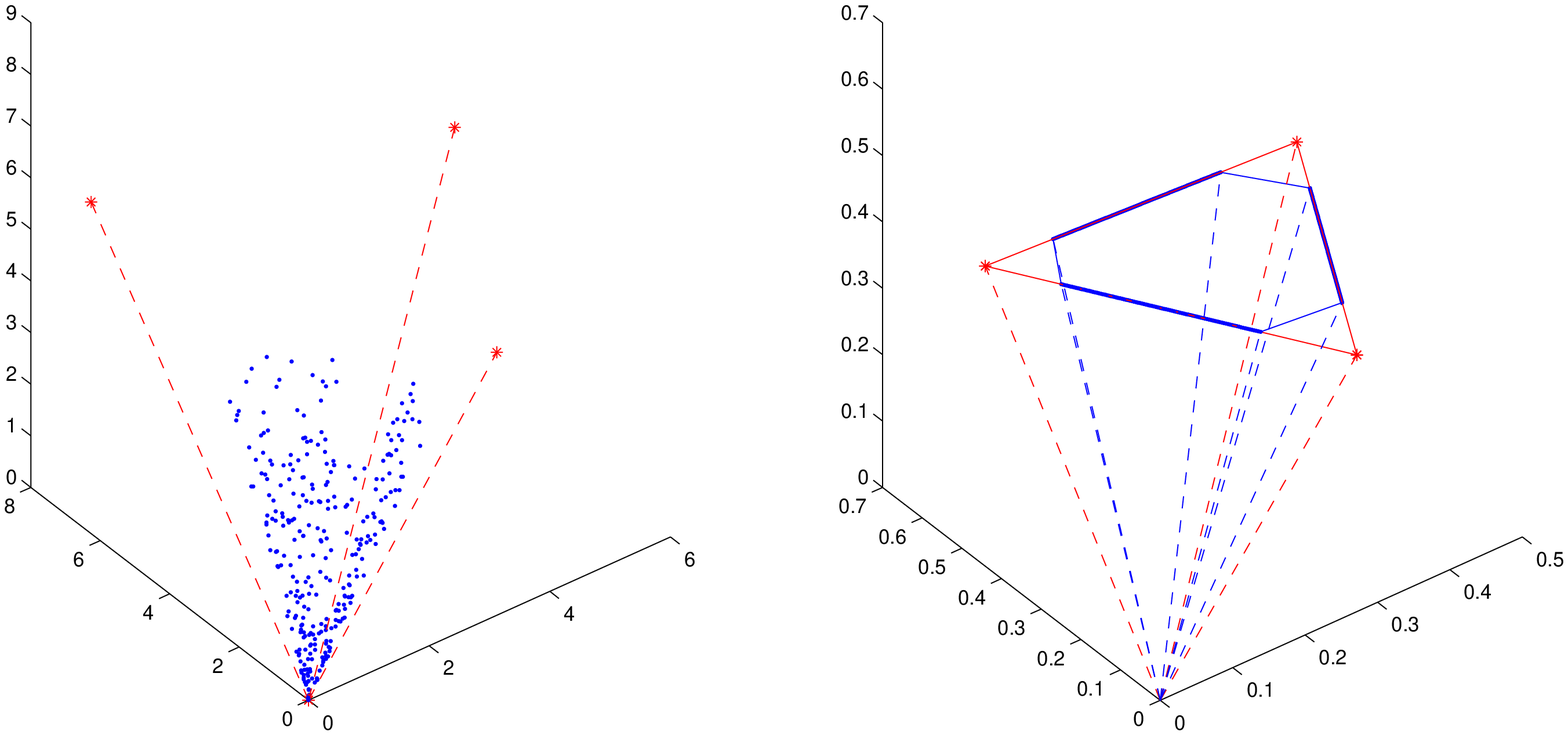}
\caption{Left: the cone $\mathcal{A}$ and the data points. Right: the convex hull $Conv(A)$ and
the convex hull $Conv(X)$ formed by rescaled data points.}
\label{cone}
\end{figure}

One of the facets is contained in the plane $\textbf{x}^{T}\cdot\textbf{1} = 1$, and
we call it \textbf{trivial facet}. The other $l-1$ facets having $\textbf{O}$ as one of
their vertices are called \textbf{nontrivial facet}. By {Theorem} \ref{thm 1},
$\mathcal{X}$ has at least the same number of nontrivial facets as $\mathcal{A}$.
In fact, it has more nontrivial facets than $\mathcal{A}$ in most cases.
If we randomly choose $m$ nontrivial facets from those of $\mathcal{X}$,
the solutions are clearly nonunique. Hence we need an additional source
assumption to provide a selection criterion:

\begin{ass}\label{ass2}
The $m$ facets of $\mathcal{X}$ containing the largest numbers of data points are contained in
the $m$ facets of $\mathcal{A}$.
\end{ass}

By assumption \ref{ass2}, we count and sort the number of data points in
each nontrivial facets of $Conv(X)$ followed by selecting the $m$ facets
with the largest numbers of data points.
Each of the facets is actually contained in one facet of $\mathcal{A}$.
So the intersection of any $m-1$ facets out of the acquired $m$ facets is
an edge of $\mathcal{A}$. By intersecting all $m$ edges
with the hyperplane $\textbf{x}^{T}\cdot\textbf{1} = 1$, we obtain the $m$ column vectors
of the mixing matrix $A$.  Finally, nonnegative least squares method yields the source matrix $S$.

Based on the ideas above, we now include the additive noise and put forward the algorithm under assumptions
\ref{nonNN}-\ref{ass2} as follows:

\begin{alg}[Face Component Analysis]\label{alg}

$(A,S)$ = FCA$(X , \rho, \epsilon , \sigma , \delta)$; parameters $\rho > 0$; $\epsilon , \sigma , \delta \in (0,1)$.
\begin{itemize}
\item[1.](Preprocessing)
$X_{0} = \max\{X,0\}$. If $\|X_{0}^{j}\|_{2} < \rho$, delete $X_{0}^{j}$ from $X_{0}$, denote by
$\hat{X_{0}}$ the resulting matrix. Project each $\hat{X_{0}}^{j}$ onto the plane $\textbf{x}^{T}\cdot\textbf{1} = 1$.

\item[2.](Convex hull)
Add the origin $\textbf{O}$ as the first column to $\hat{X_{0}}$.
Return all the facets and vertices of $Conv(X)$, keep only the nontrivial facets. Denote by $F_{i}$
the $i$th facet and $V_{i}$ the set of its vertices.

\item[3.](Grouping)
Initialize $G_{i} = V_{i}$. If $\hat{X_{0}}^{j} \not\in G_{i}$ and $\textrm{dist}(\hat{X_{0}}^{j},F_{i}) < \epsilon \approx 0$ and $ \textrm{dist} (\hat{X_{0}}^{j},V_{i}) > \sigma \approx 0$, add $\hat{X_{0}}^{j}$ to $G_{i}$.

\item[4.](Plane fitting)
For each $G_{i}$, obtain the equation of its fitting plane denoted as $\textbf{x}^{T}\cdot \textbf{b}_{i} = 0$, where $\textbf{b}_{i}$ is the normal vector with norm $1$. Select $m$ planes from $G_{i}$'s with the largest cardinalities such that the inner product of any two $\textbf{b}_{i}$s $< \delta \approx 1$.

\item[5.](Intersecting)
Obtain the $m$ intersections of any $m-1$ planes out of $m$ planes from step 4 with
the plane $\textbf{x}^{T}\cdot\textbf{1} = 1$, and form the $m$ columns of mixing matrix $A$.

\item[6.](Source recovering)
For each column $X_{0}^{j}$, solve the following optimization problem to find
the corresponding column $S^{j}$ of $S$ as:

\textbf{minimize} \quad $\|X_{0}^{j} - A\,S^{j}\|_{2}$,

\textbf{subject to } \quad $S^{j} \succcurlyeq 0$
\end{itemize}
\end{alg}

\begin{figure}
\includegraphics[height=7cm,width=15cm]{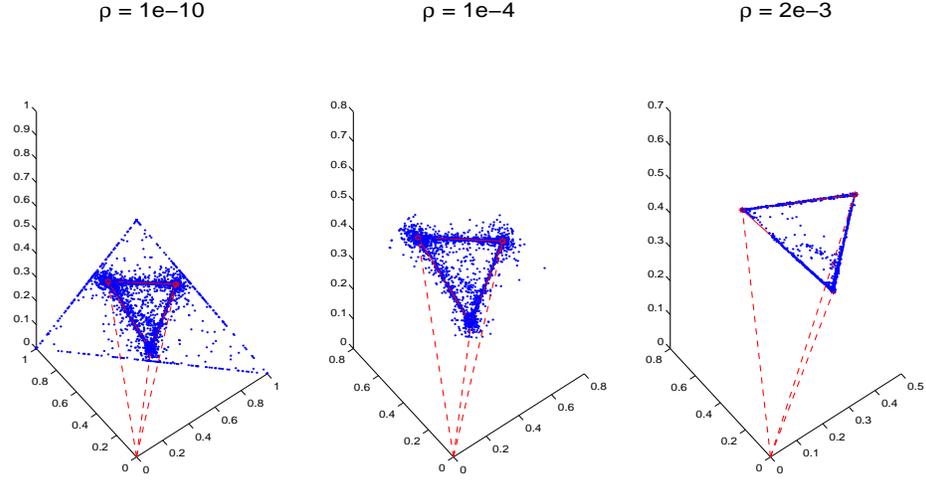}
\caption{Comparison of denoising effects with different $\rho$.}
\label{denoise}
\end{figure}

\begin{rem}

\begin{itemize}

\item[1.]
We can denoise by increasing the threshold $\rho$ in step 1. The effects of denoising with varied $\rho$ are present in Fig. \ref{denoise}. We add white Gaussian noise with signal-to-noise ratio (SNR) = $80$ dB to three noiseless mixtures. In the left plot of Fig. \ref{denoise}, an outside triangle forms due to the intersections of the plane $x+y+z = 1$, x-plane, y-plane and z-plane. When noiseless signal is corrupted by additive white Gaussian noise, many data points with norms of almost $0$ tend to get negative entries.  Then we apply the threshold of $0$ to $X$'s entries in step 1. Thus these points fall into x-plane, y-plane or z-plane. If we choose a very very small $\rho = 10^{-10}$ as in the left plot, these points will still remain (i.e. hardly any denoising effect). So projecting them onto the plane $x+y+z = 1$ forms the outside triangle. It is awful to get such a geometric structure because it is impossible to identify $\mathcal{A}$.
However, if $\rho$ goes up to $2\times10^{-3}$ as in the right plot, the structure of $\mathcal{A}$ emerges. Clearly a better geometric structure of data points is achieved by larger $\rho$. However, we also need to avoid $\rho$ being too large in that we may lose the structure of $Conv(X)$ with few data points.

\item[2.]
Considering the presence of noise, we modify the criterion of one point belonging to a facet as follows.
The point is not required to lie exactly on the facet.
Instead, we only require that it is near the facet yet
not around the vertices of the facet. The thresholds $\epsilon$ and $\sigma$ in the following
section of numerical experiments vary from $10^{-6}$ to $0.1$ depending on the level of noises. If there is almost no noise, we can achieve perfect recovery by setting these two thresholds to be any values less than $10^{-5}$.
A higher level of noises demands larger values of $\epsilon$ and $\sigma$.
These two parameters are usually of the same order.

\item[3.]
In step 4, we introduce the threshold $\delta$ to avoid selecting two nearly
coplanar fitting planes which may actually correspond to the same facet of $Conv(X)$.
Normally we set $\delta$ to be $0.99$.

\item[4.] We will introduce denoising methods by smoothing filters such as box filter, Gaussian filter, and total variation (TV denoising) in next section.  The combination of these denoising methods and our FCA tend to perform better than FCA alone when the noise level is high.  We embed the noise filter into FCA after step 3 where we are able to preserve the data points close to the facets of $Conv(X)$ only.  Then applying the denoising methods yield a more desirable geometric structure of data points.
\end{itemize}
\end{rem}

\section{Numerical Experiments}
We report the numerical results of our algorithm in this section.  The data we have tested
include real-world NMR spectra as well as synthetic mixtures.
All the entries of the mixing matrices and the sources matrices are positive.

\begin{figure}
\includegraphics[height=7cm,width=15cm]{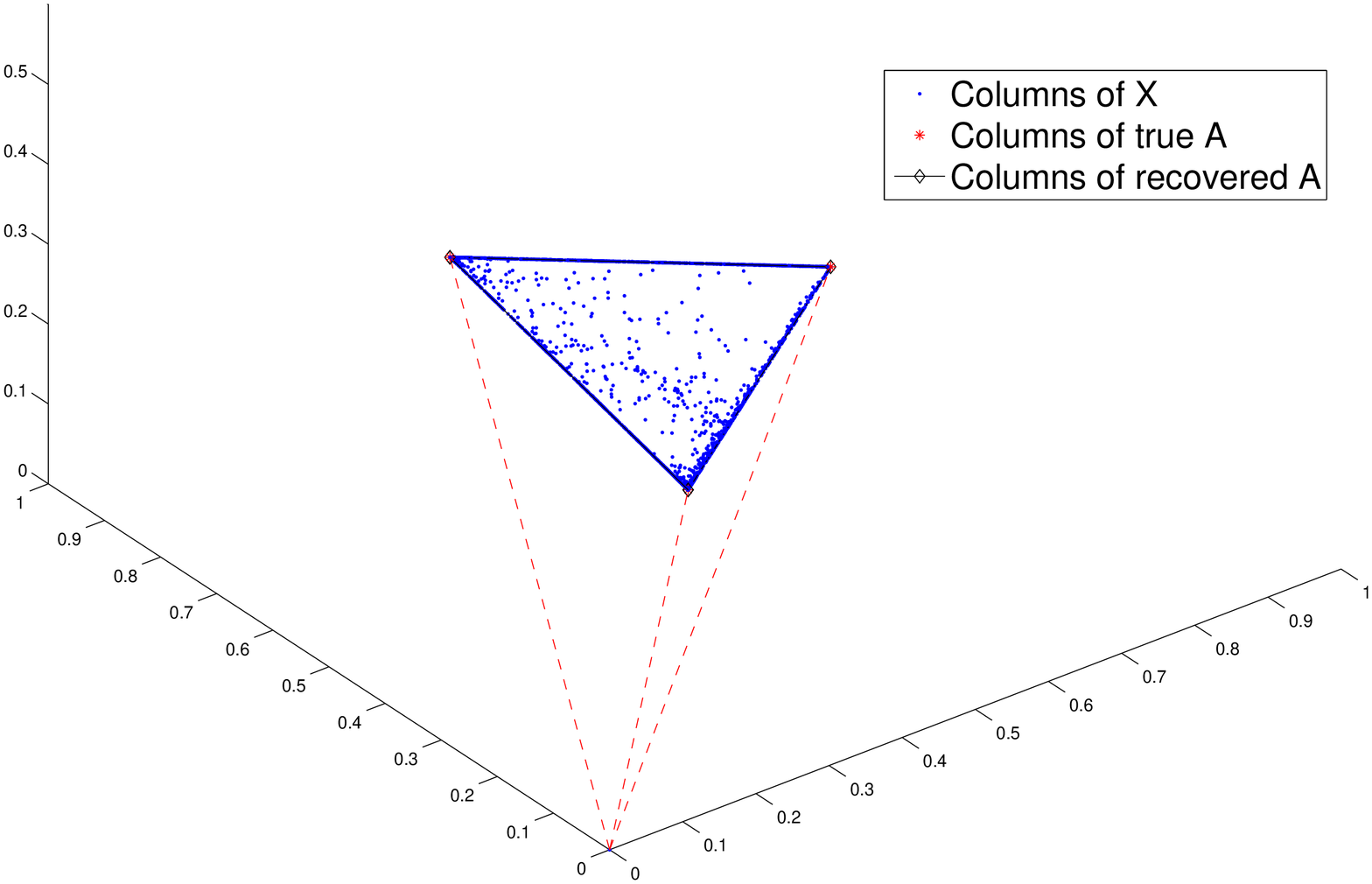}
\caption{Rescaled columns of mixture matrix and columns of mixing matrix from Example 1.}
\label{NNpoint}
\end{figure}

\begin{figure}
\includegraphics[height=7cm,width=15cm]{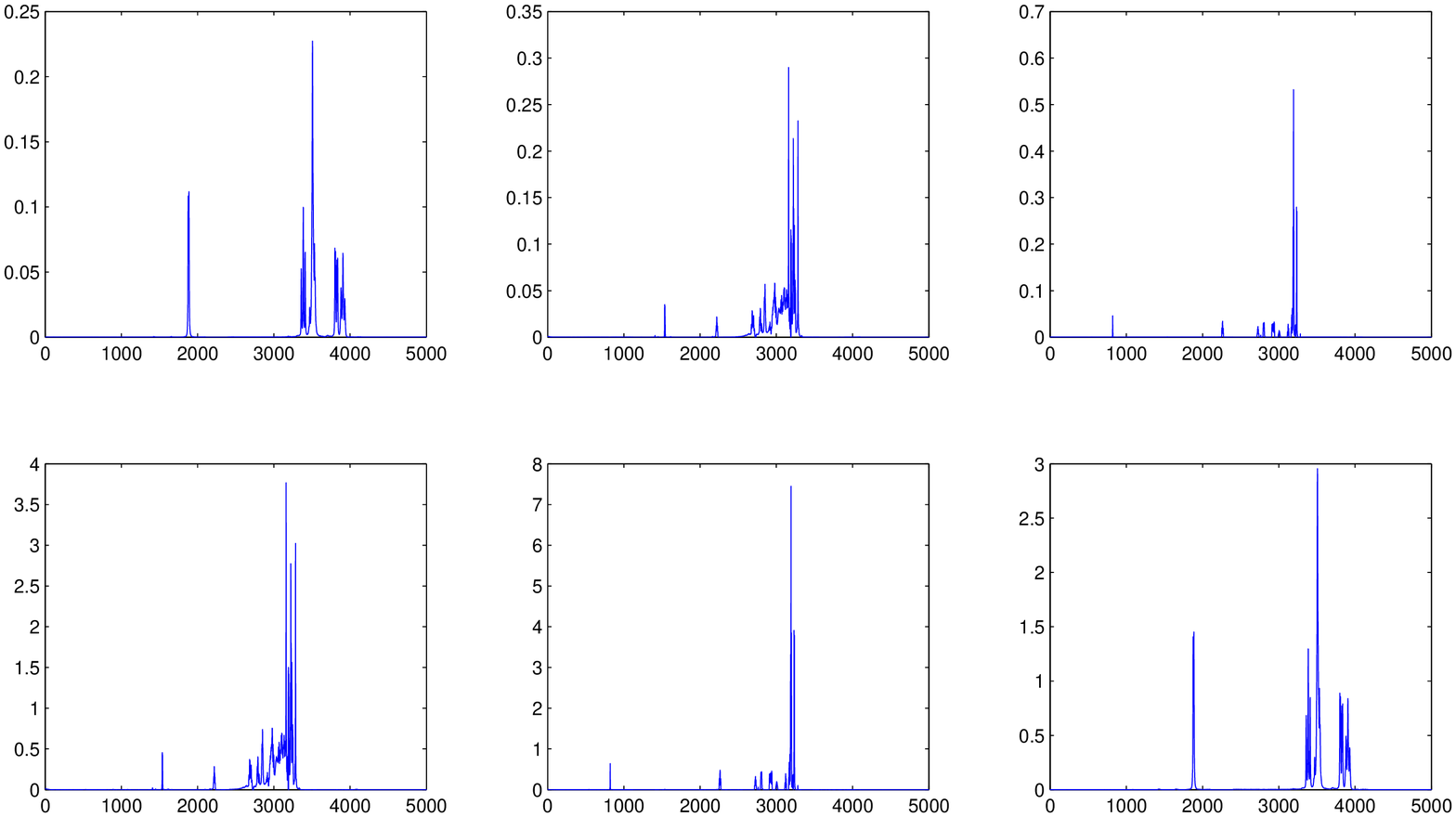}
\caption{Top row: from left to right, the three reference spectra of $\beta$-cyclodextrine, $\beta$-sitosterol, and menthol.  Bottom row: recovery results by our method.}
\label{NNsources}
\end{figure}

As we have pointed out in the last section,
NN's assumption is a special case of ours. So our algorithm is supposed
to work for the separations of NN data.
The first example is to recover three sources from three mixtures,
where the source signals have stand-alone peaks.
For the data, we used true NMR spectra of four compounds
$\beta$-cyclodextrine,$\beta$-sitosterol, and menthol as source signals.  The NMR spectrum of a chemical compound is produced by the Fourier transformation of a time-domain
signal which is a sum of sine functions with exponentially decaying envelopes \cite{Ern}.  The real part of the spectrum can be presented as the sum of symmetrical, positive valued, Lorentzian-shaped peaks.  The NMR reference spectra of $\beta$-cyclodextrine,$\beta$-sitosterol, and menthol are shown in the top panel of Fig. \ref{NNsources} from left to right.  For the parameters, we set $\rho = 10^{-3}, \epsilon = \sigma = 10^{-6}, \delta = 0.99$.  Fig. \ref{NNpoint} shows the geometric structure of the mixtures and mixing matrix.  The reference spectra and computational results are shown in Fig. \ref{NNsources}. $A_{1}$ is the rescaled true mixing matrix, while $\hat{A_{1}}$ is the computed mixing matrix via our method. Apparently the separation results are very nice in that $A_{1}$ and $\hat{A_{1}}$ are identical.

$$
A_{1} =
\left(
  \begin{array}{ccc}
    0.0769  &  0.4615  &  0.3571 \\
    0.3846  &  0.4615  &  0.0714 \\
    0.5385  &  0.0769  &  0.5714 \\
  \end{array}
\right)
$$

$$
\hat{A_{1}} =
\left(
  \begin{array}{ccc}
    0.4615  &  0.3571  &  0.0769 \\
    0.4615  &  0.0714  &  0.3846 \\
    0.0769  &  0.5714  &  0.5385 \\
  \end{array}
\right)
$$

\begin{figure}
\includegraphics[height=7cm,width=15cm]{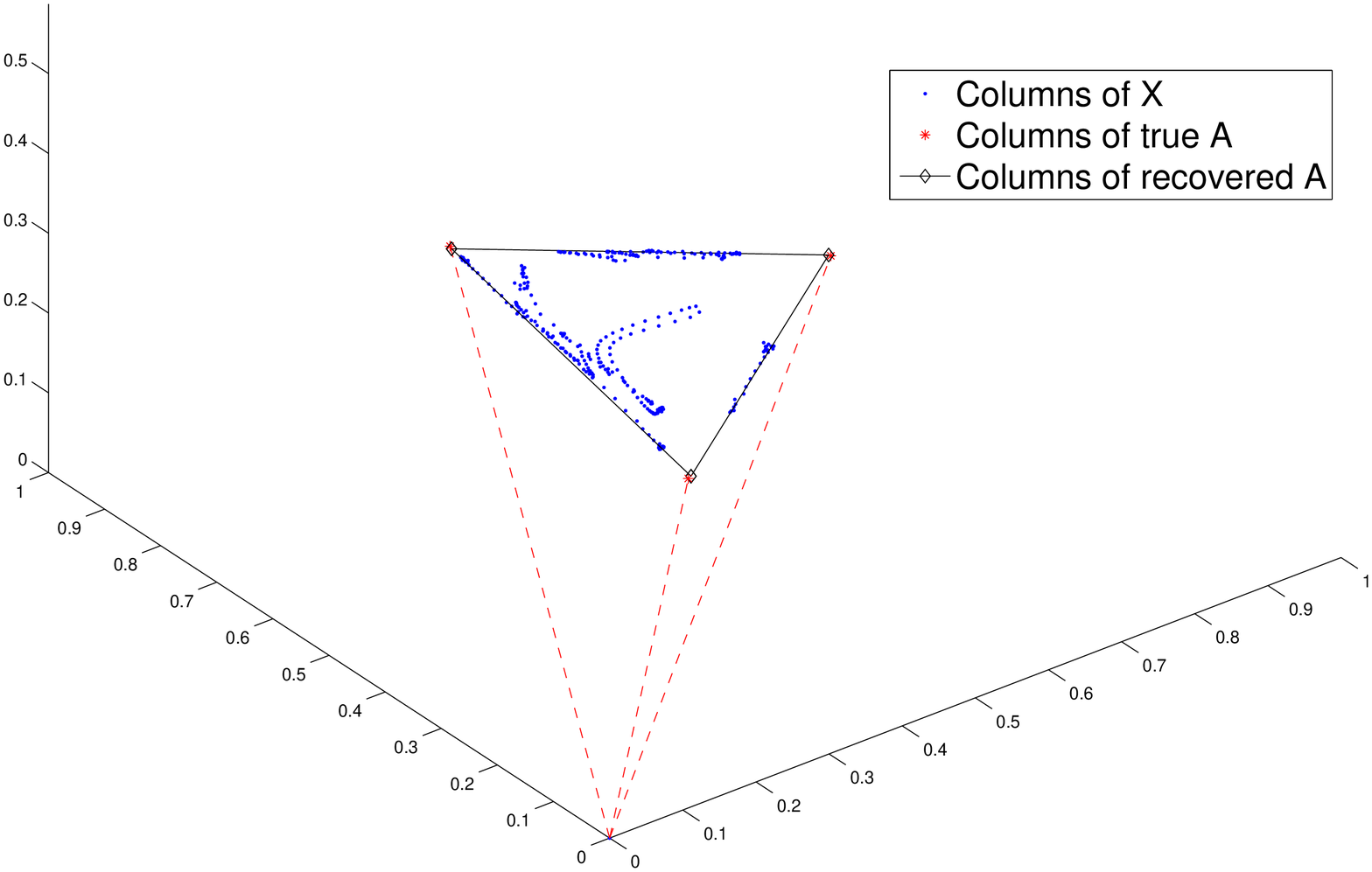}
\caption{Rescaled columns of mixture matrix and columns of mixing matrix from Example 2.
Some noisy data points are already deleted from $X$ by step 1 of Algorithm \ref{alg}.
Parameters are $\rho = 50, \epsilon = 5\times10^{-3}, \sigma = 6\times10^{-3}, \delta = 0.99$.}
\label{nonNNpoint}
\end{figure}

\begin{figure}
\includegraphics[height=7cm,width=15cm]{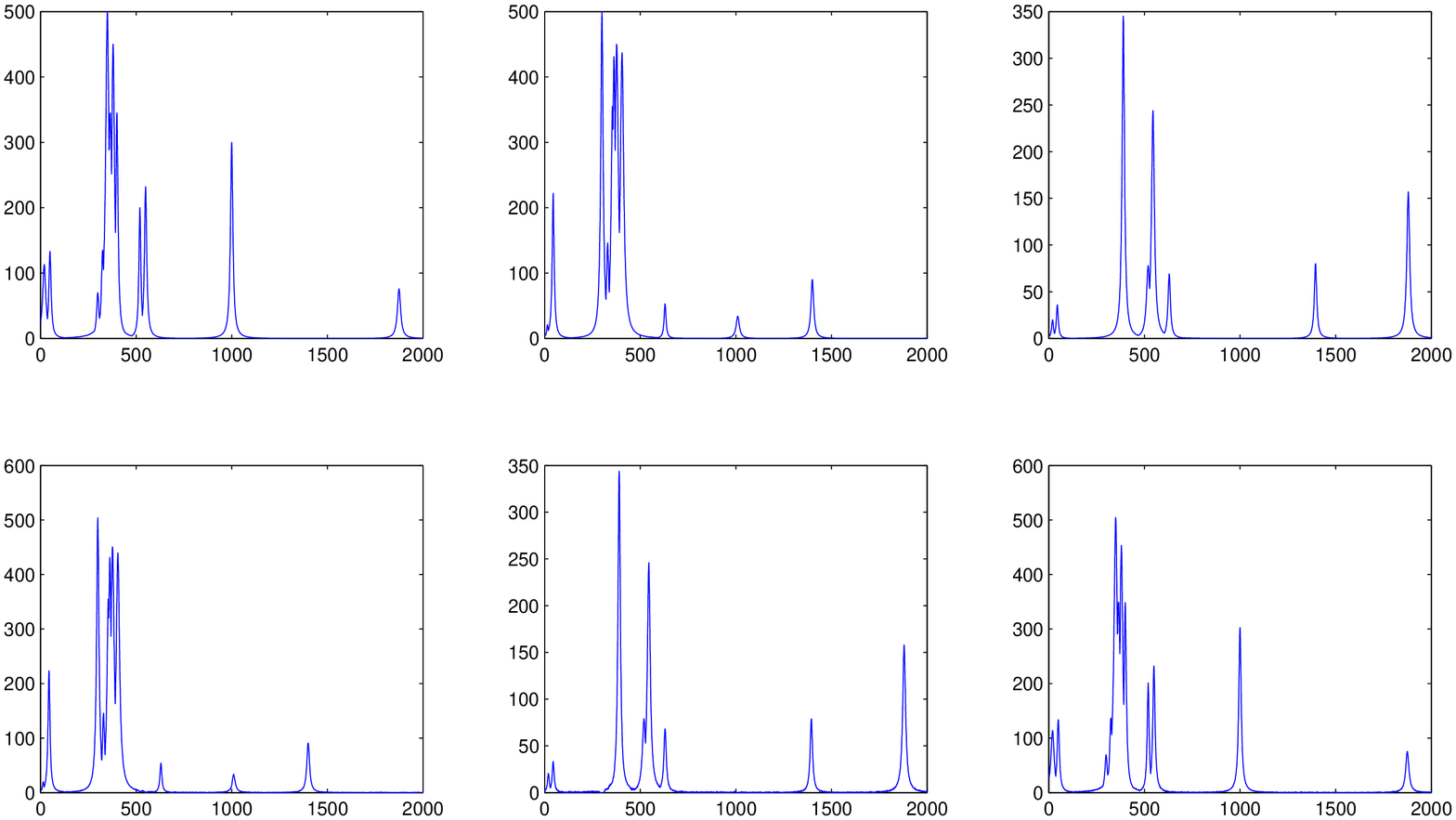}
\caption{Top row: the true source signals. Bottom row: computed source signals via our method. Prameters are $\rho = 50, \epsilon = 5\times10^{-3}, \sigma = 6\times10^{-3}, \delta = 0.99$.}
\label{nonNNsources}
\end{figure}

In a second example, there are three Lorentzian source signals which no longer satisfy NN assumption.
We use them to create three noisy mixtures by adding white Gaussian noise with SNR = 50 dB.
A good result is achieved by setting $\rho = 50, \epsilon = 5\times10^{-3}, \sigma = 6\times10^{-3}, \delta = 0.99$.
Fig. \ref{nonNNpoint} and Fig. \ref{nonNNsources} show the geometry of data points and
recovery results. True mixing matrix $A_{2}$ and computed $\hat{A_{2}}$ are as below (first row of $\hat{A_2}$ is scaled to be same as that of $A_2$)

$$
A_{2} =
\left(
  \begin{array}{ccc}
    0.0769  &  0.4615  &  0.3571 \\
    0.3846  &  0.4615  &  0.0714 \\
    0.5385  &  0.0769  &  0.5714 \\
  \end{array}
\right),
$$

$$
\hat{A_{2}} =
\left(
  \begin{array}{ccc}
    0.4615  &  0.3571  &  0.0769 \\
    0.4565  &  0.0729  &  0.3700 \\
    0.0823  &  0.5765  &  0.5106 \\
  \end{array}
\right).
$$

\begin{figure}
\includegraphics[height=7cm,width=15cm]{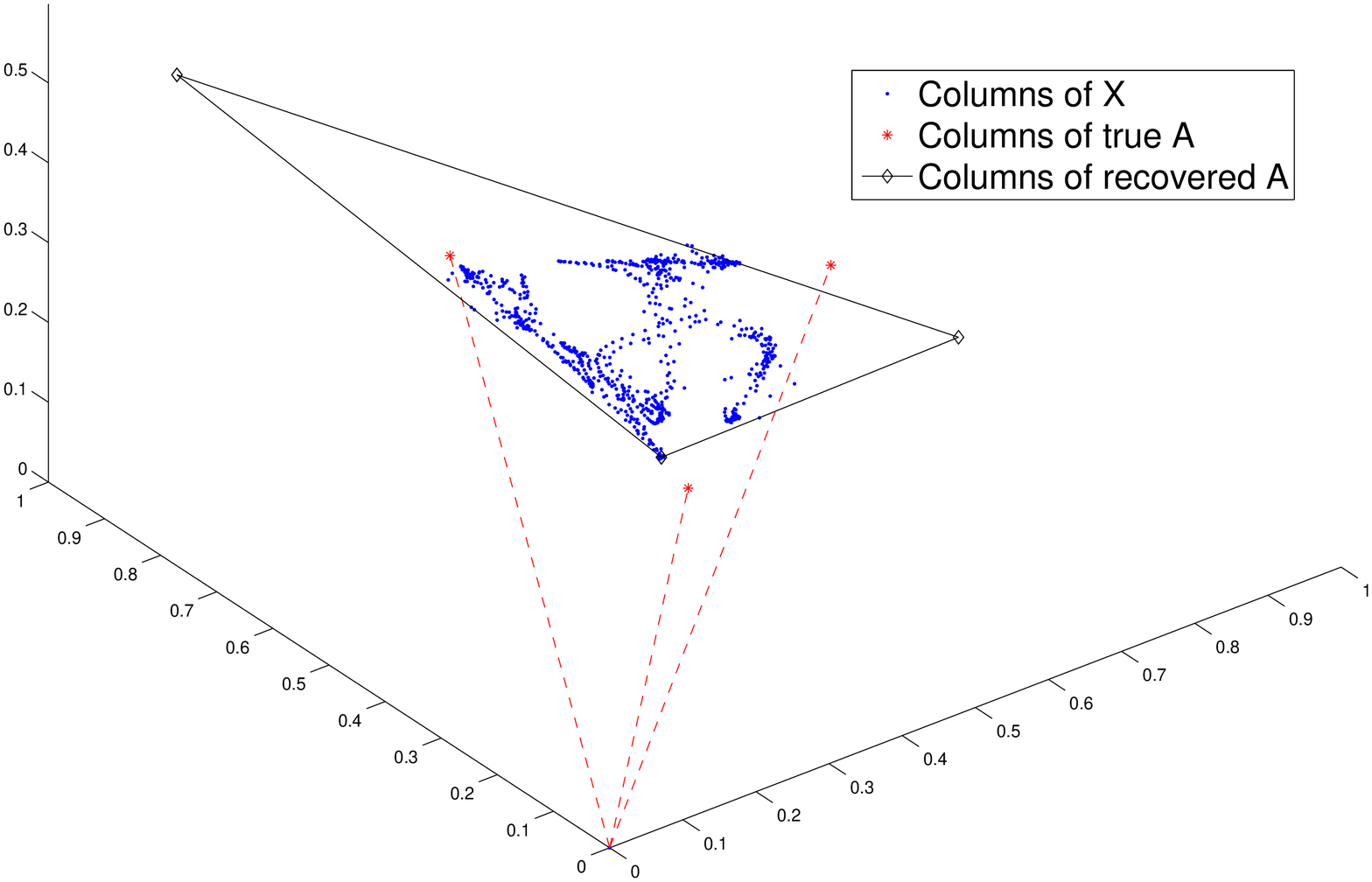}
\caption{Example 2: $\rho = 5, \epsilon = 5\times10^{-3}, \sigma = 6\times10^{-3}, \delta = 0.99$.}
\label{nonNNpoint1}
\end{figure}

\begin{figure}
\includegraphics[height=7cm,width=15cm]{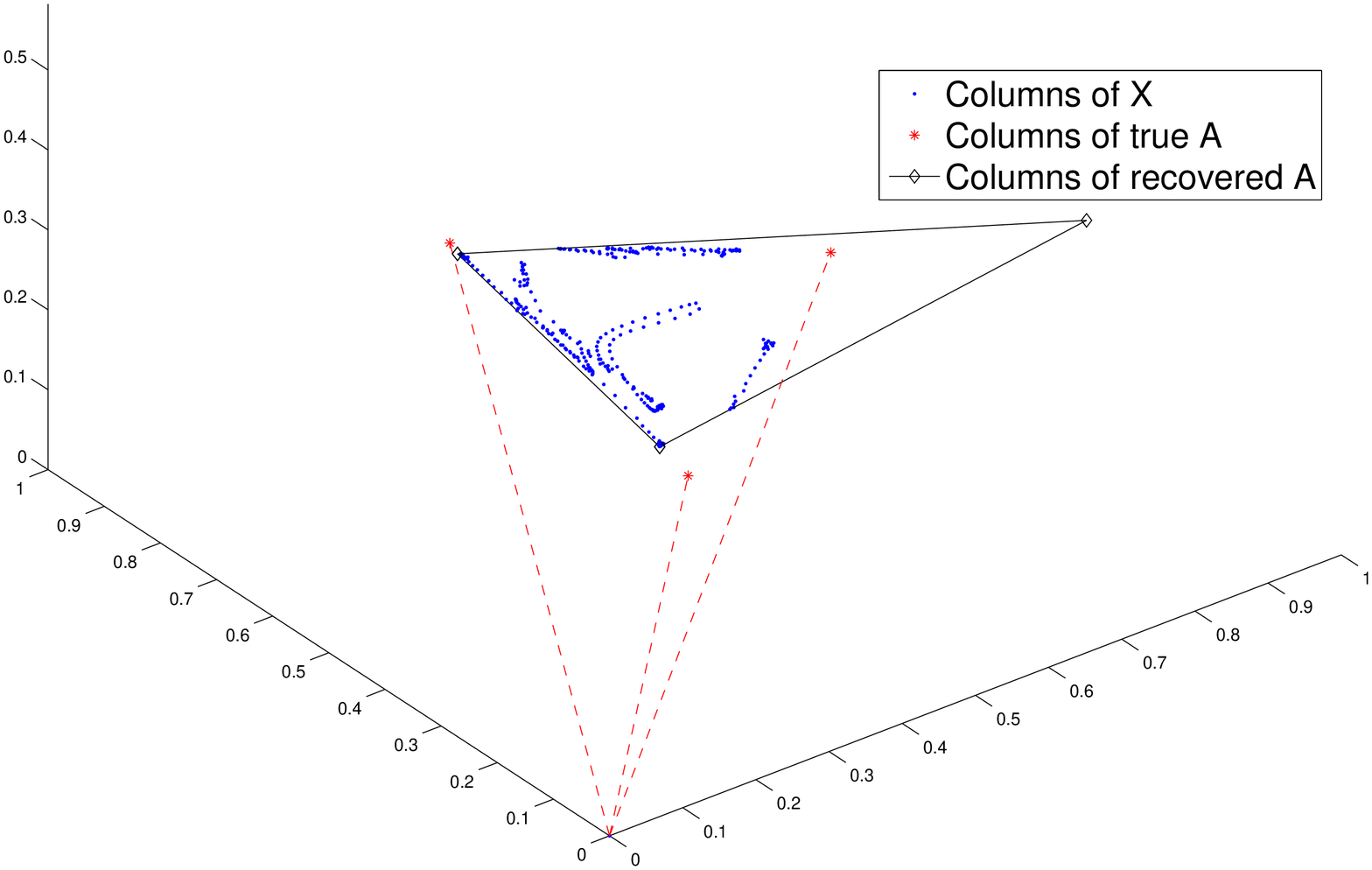}
\caption{Example 2: $\rho = 50, \epsilon = 10^{-5}, \sigma = 10^{-5}, \delta = 0.99$.}
\label{nonNNpoint2}
\end{figure}

\begin{figure}
\includegraphics[height=7cm,width=15cm]{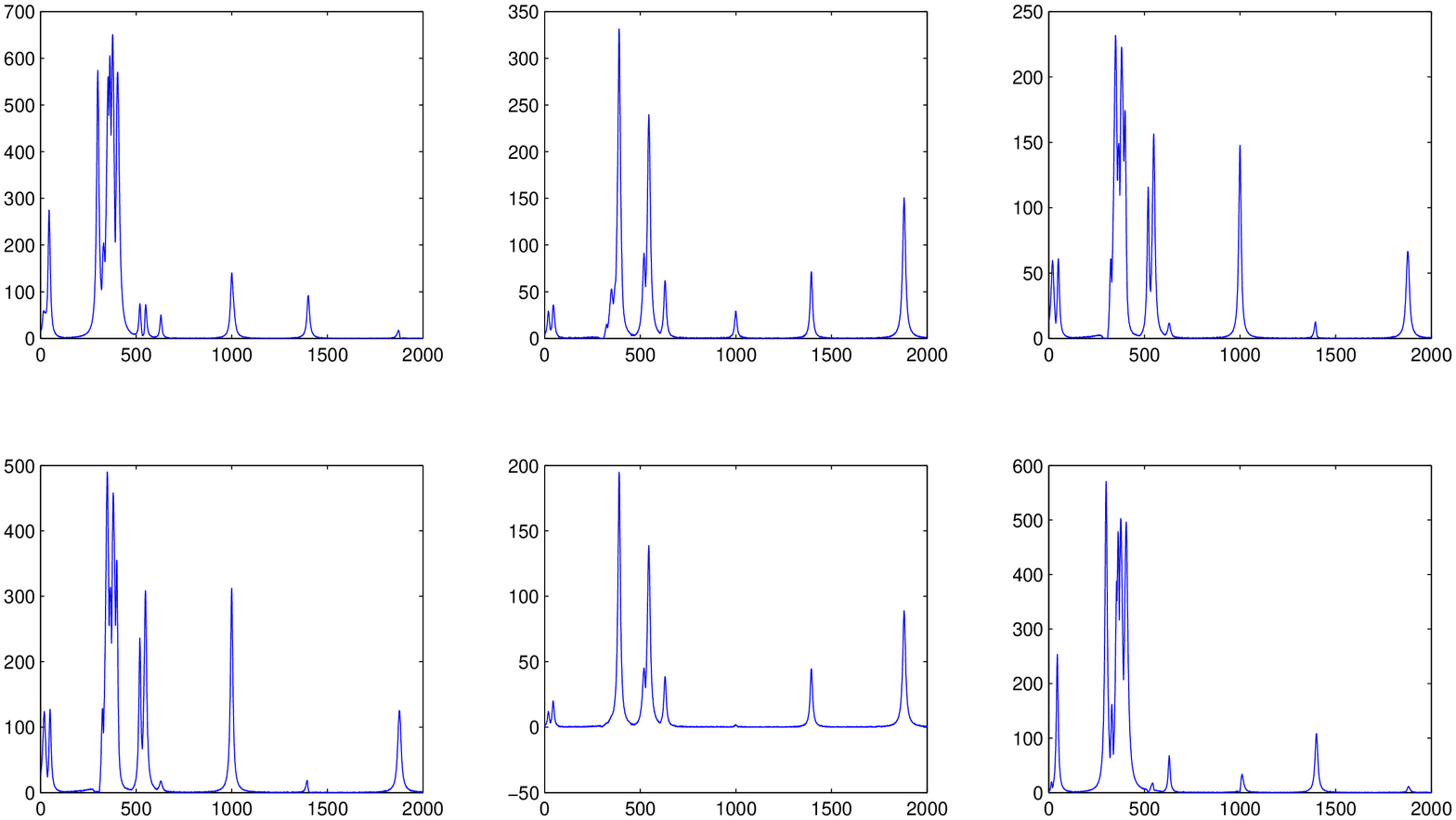}
\caption{Computed sources: Top row: $\rho = 5, \epsilon = 5\times10^{-3}, \sigma = 6\times10^{-3}, \delta = 0.99$. Bottom row:  $\rho = 50, \epsilon = 10^{-5}, \sigma = 10^{-5}, \delta = 0.99$.}
\label{nonNNsources1}
\end{figure}

To provide further insight into how to choose the parameters, two results caused by
inappropriate threshold values are presented as well in Fig. \ref{nonNNpoint1}, Fig. \ref{nonNNpoint2}
and Fig. \ref{nonNNsources1}. In the first experiment, some noisy data points destroy the
geometric structure as shown in Fig. \ref{nonNNpoint1} since $\rho$ is not big
enough to filter out them. In the second one, $\epsilon$ and $\sigma$ are
too small for the level of noises.

\begin{figure}
\includegraphics[height=7cm,width=15cm]{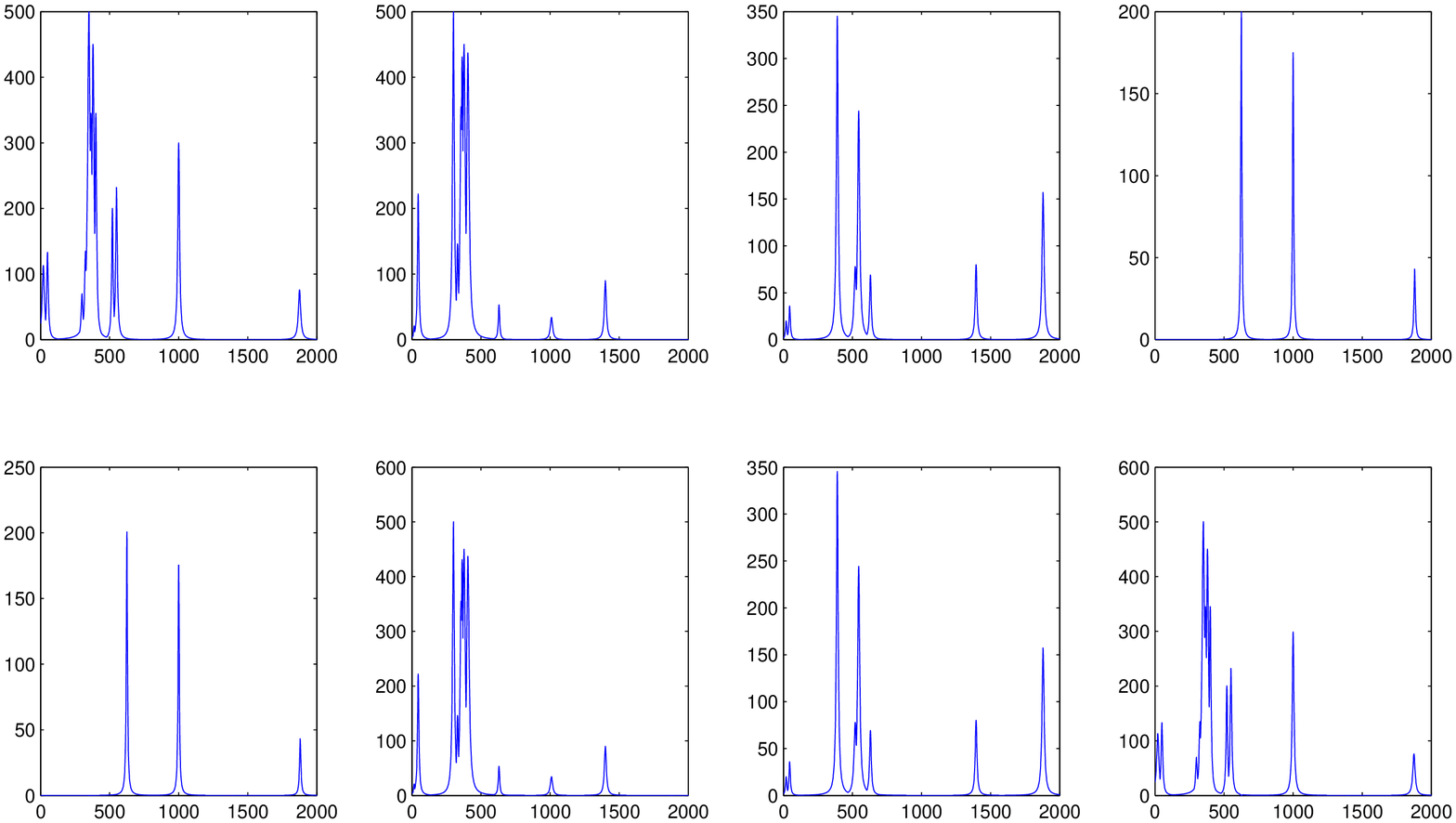}
\caption{Top row: the four original source signals. Bottom row: the recovered sources.}
\label{Sources}
\end{figure}

The last example concerns BSS in higher dimensions.
We manage to recover four sources from four mixtures.
The data points are in 4-dimensional space, so they are difficult to visualize.
We chose $\rho = 1, \epsilon = 2\times10^{-5}, \sigma = 10^{-5}, \delta = 0.99$ here.
The original sources and recovered ones are present in Fig. \ref{Sources}.
The true and computed mixing matrices are shown below (first row of $\hat{A_3}$ is scaled to be same as that of $A_3$).
$$
A_{3} =
\left(
  \begin{array}{cccc}
    0.1923 & 0.2500 & 0.2632 & 0.1000 \\
    0.1923 & 0.2500 & 0.2105 & 0.2000 \\
    0.2692 & 0.3750 & 0.4211 & 0.3000 \\
    0.3462 & 0.1250 & 0.1053 & 0.4000 \\
  \end{array}
\right)
$$

$$
\hat{A_{3}} =
\left(
  \begin{array}{cccc}
    0.1000 & 0.2500 & 0.2632 & 0.1923 \\
    0.1997 & 0.2500 & 0.2107 & 0.1922 \\
    0.2992 & 0.3749 & 0.4211 & 0.2694 \\
    0.4011 & 0.1252 & 0.1057 & 0.3456 \\
  \end{array}
\right)
$$

\begin{figure}
\includegraphics[height=6cm,width=15cm]{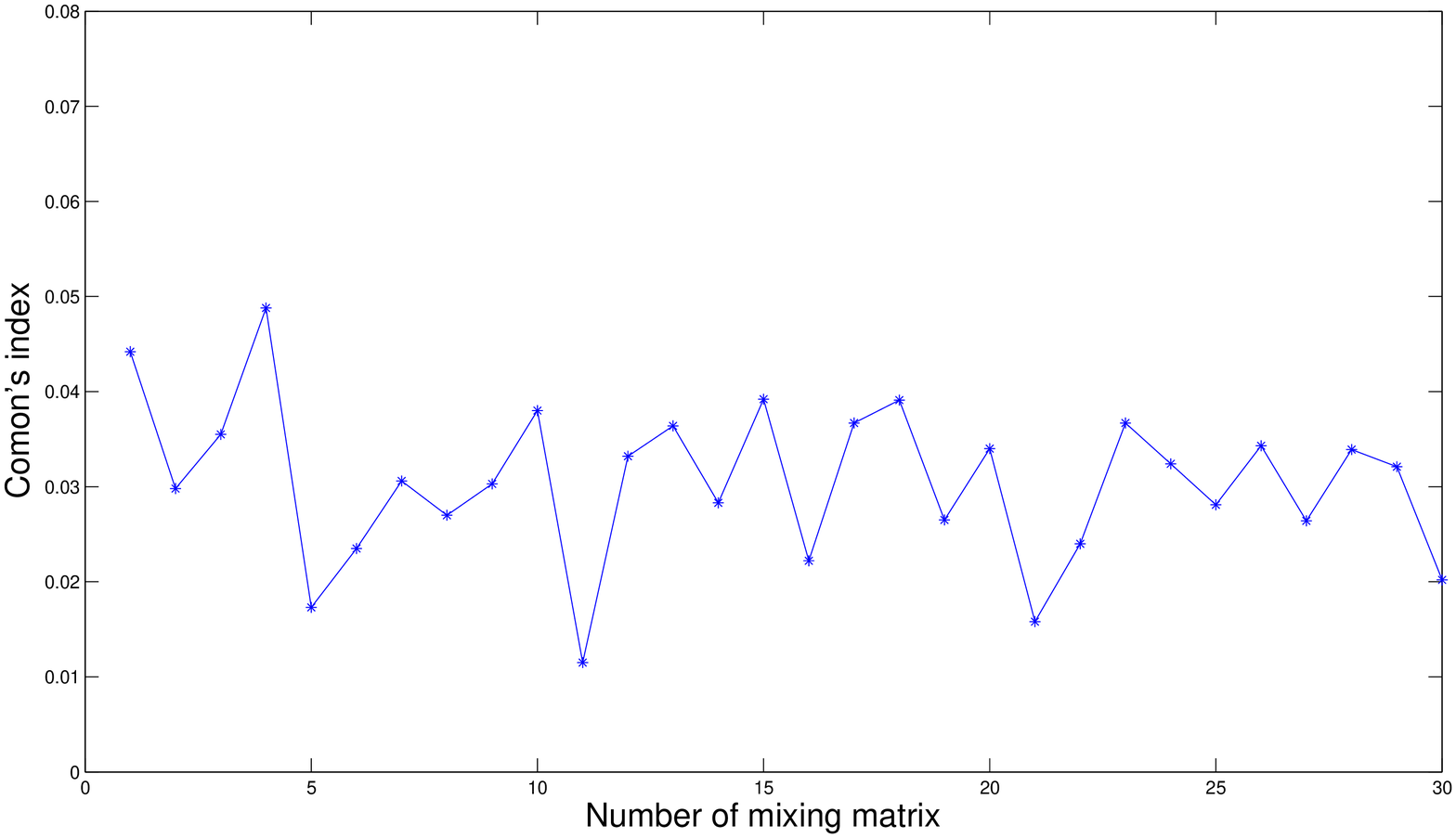}
\caption{Performance of our method on 30 random $4\times 4$ mixing matrices. The four sources in Example 3 are used.}
\label{idx1}
\end{figure}

\begin{figure}
\includegraphics[height=6cm,width=15cm]{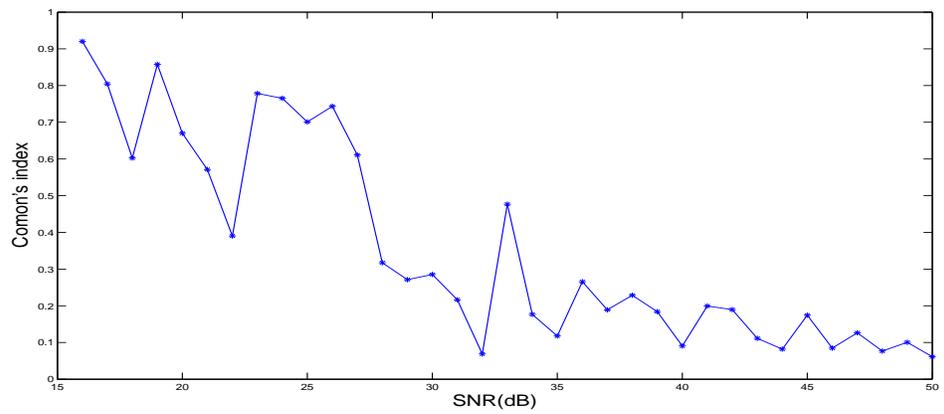}
\caption{Test the robustness of our method in the presence of noise.}
\label{idx2}
\end{figure}

To test the performance of our method, we compute the Comon's index \cite{Comon}.
The index is defined as follows:
let $A$ and $\hat{A}$ be two
nonsingular matrices with $L_{2}$-normalized columns.  Then the distance between $A$ and $\hat{A}$
denoted by $\epsilon(A,\hat{A})$ which reads
\begin{equation*}
\epsilon(A,\hat{A}) = \sum_i\biggl | \sum_j |d_{ij}|-1\biggr |^2 + \sum_j\biggl | \sum_i |d_{ij}|-1\biggr |^2
 + \sum_i\biggl | \sum_j |d_{ij}|^2-1\biggr | + \sum_j\biggl | \sum_i |d_{ij}|^2-1\biggr |\;,
\end{equation*}
where $D =A^{-1}\hat{A} $, and $d_{ij}$ is the entry of $D$.  In \cite{Comon} Comon proved that
$A$ and $\hat{A}$ are considered nearly equivalent in the sense of BSS (i.e., $\hat{A} = A\,P\,\Lambda$)
if $\epsilon(A,\hat{A}) \approx 0$. Fig. \ref{idx1} and Fig. \ref{idx2} show Comon's indices
between the true mixing matrices and the computed matrices by our method.
For the result in Fig. \ref{idx1}, we compute the Comon's indices using the four sources
in Example 3 and 30 $4\times 4$ random mixing matrices.
Clearly the Comon's indices are very small suggesting the equivalence in the sense of
BBS of the true mixing matrices and the computed ones. Fig. \ref{idx2} shows the performance of
our method in the presence of noise. The three sources in Example 2 are
combined to generate three noisy mixtures by adding white Gaussian noises
with SNR varying from 16 dB to 50 dB. The reliability of our method is manifested from the plot.

\subsection{Denoising}
If there is considerable noise in the data, it would be desirable to reduce or remove the noise before feeding them to the proposed method.  We shall propose to apply imaging denoising techniques to the FCA for noise reduction or removal. These denoising methods may be used when the noise level is high, and they can be combined with the FCA after step 3.
In image processing, the simplest denoising method is the {\em sliding mean} or {\em box filter} \cite{MCD}.  Each pixel value is replaced by the mean of its local neighbors.  The Gaussian filter is similar to the box filter, except that the values of the neighboring pixels are given different weighting, that being defined by a spatial Gaussian distribution.  The Gaussian filter is probably the most widely used noise reducing filter.

\begin{figure}
\includegraphics[height=7cm,width=7.5cm]{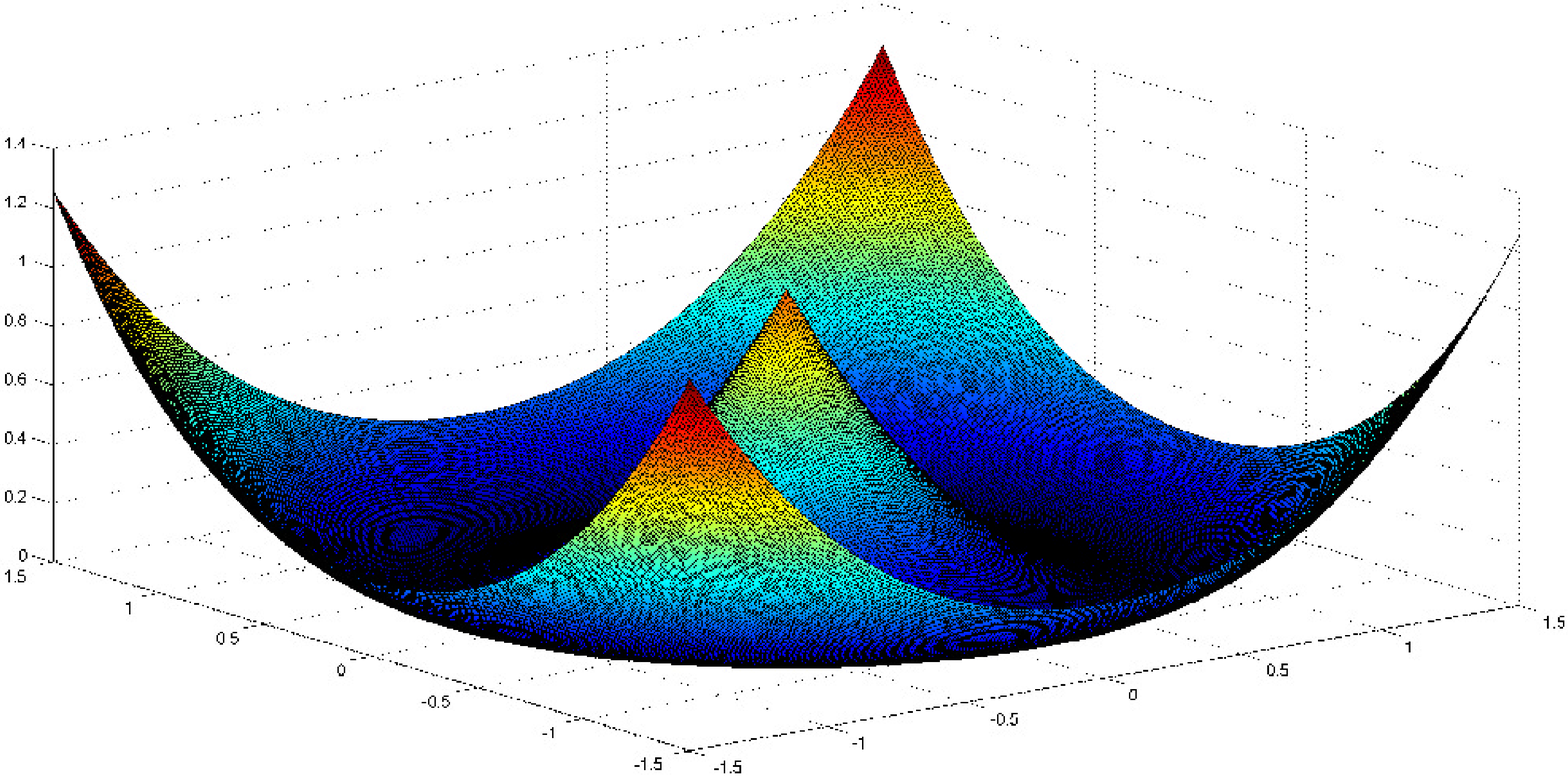}
\includegraphics[height=7cm,width=7.5cm]{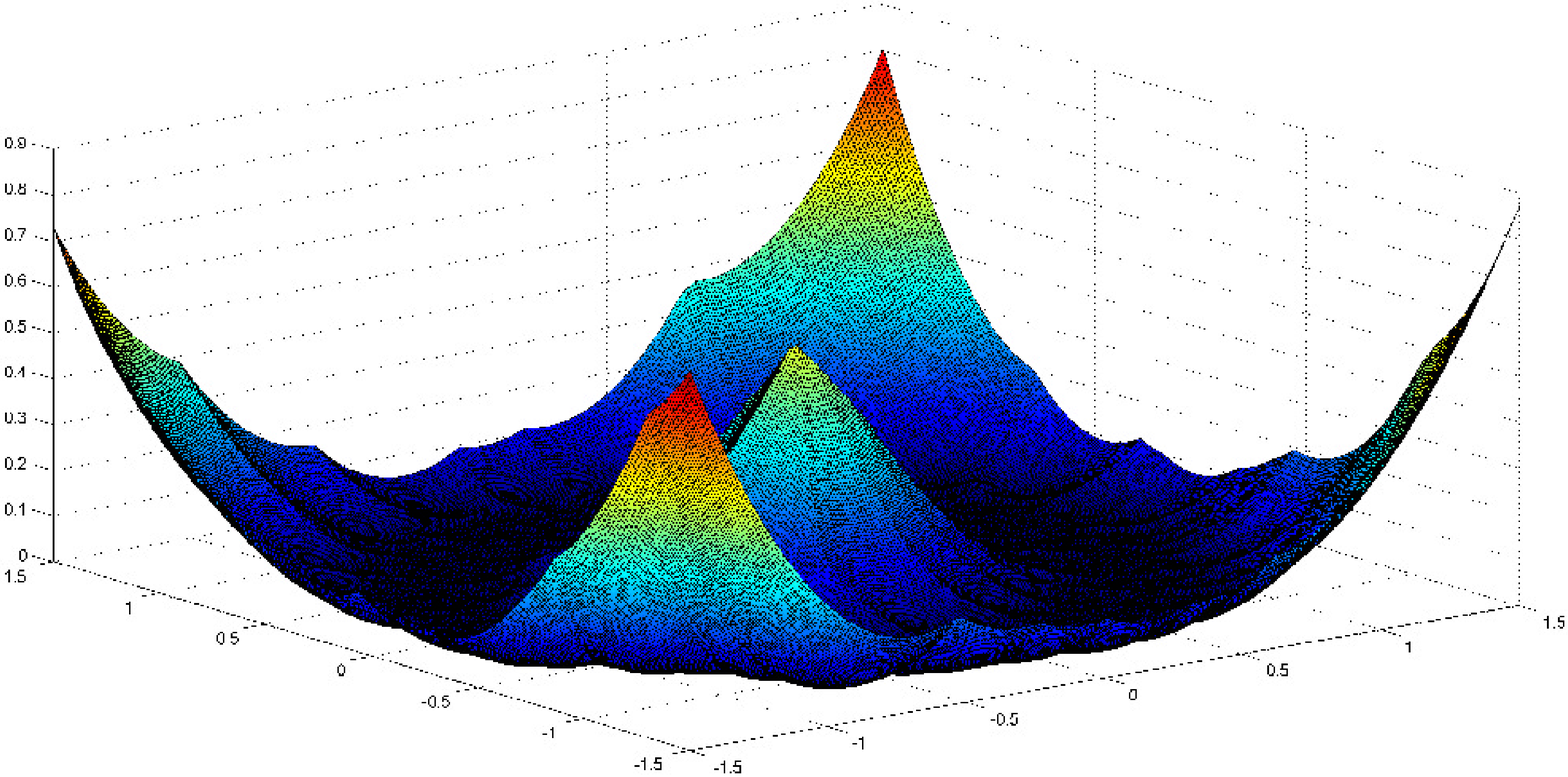}
\caption{The distance function with and without noise.}
\label{levelset}
\end{figure}

These local image smooth filters could be applied for the point cloud noise reduction with slight modification.  They are incorporated in to FCA after step 3 where the groups $G_i$ have been obtained.  Within each group, the $K$-nearest neighbors are searched based on Euclidean distance, the best choice of $K$ depends upon the data and noise level.  For example, fewer neighbors should be selected when less noise presents.  Note that these denoising methods generalize to any dimensional point cloud, they are easy to implement. A disadvantage is that they tend to smooth away edge or corner structures while reducing the noise.  To overcome this shortcoming, we propose to apply the total variation idea of image denoising for noise removal.  It is based on the principle that signals with excessive noise have high total variation, that is, the integral of the absolute gradient of the signal is high.  According to this principle, reducing the total variation of the signal subject to it being a close to the original signal, removes the unwanted detail whilst preserving important details such as edges.  The concept was originated in Rudin, Osher, and Fatemi in 1992 \cite{ROF_92}, and it can be applied to the point cloud noise removal.  In the following, we shall use the example of a point cloud in $xyz$ plane to
illustrate the idea of total variation denoising.  We first preprocess the data by rescaling
them onto a plane $x+y+z = 1$, the projected data are two dimensional.
For each point $(x,y)$, a distance function to the point cloud (data points) is defined as
 \begin{equation*}
d(x,y) = \min_{x_i,y_i}\sqrt{(x-x_i)^2 + (y-y_i)^2}\;,
\end{equation*}
where $(x_i,y_i)^{\mathrm{T}}$ corresponds to the $i$-th column of $X$, note that $z_i$ is not included since the $z_i = 1-x_i -y_i$.  Fig. \ref{levelset} shows an example of distance function to a unit circle, the left plot is the distance function with no noise, while the right plot is the function with noise.
For computation, the distance function will be restricted on a rectangular region which contains the point cloud.
 Note that $d(x,y) \geq 0$,  the distance function actually defines intensities of an image.
Clearly, the set of $\mathcal{{S}} = \{(x,y): d(x,y) = 0\}$ is the point cloud for the noiseless case.

\begin{figure}
\includegraphics[height=7cm,width=7.5cm]{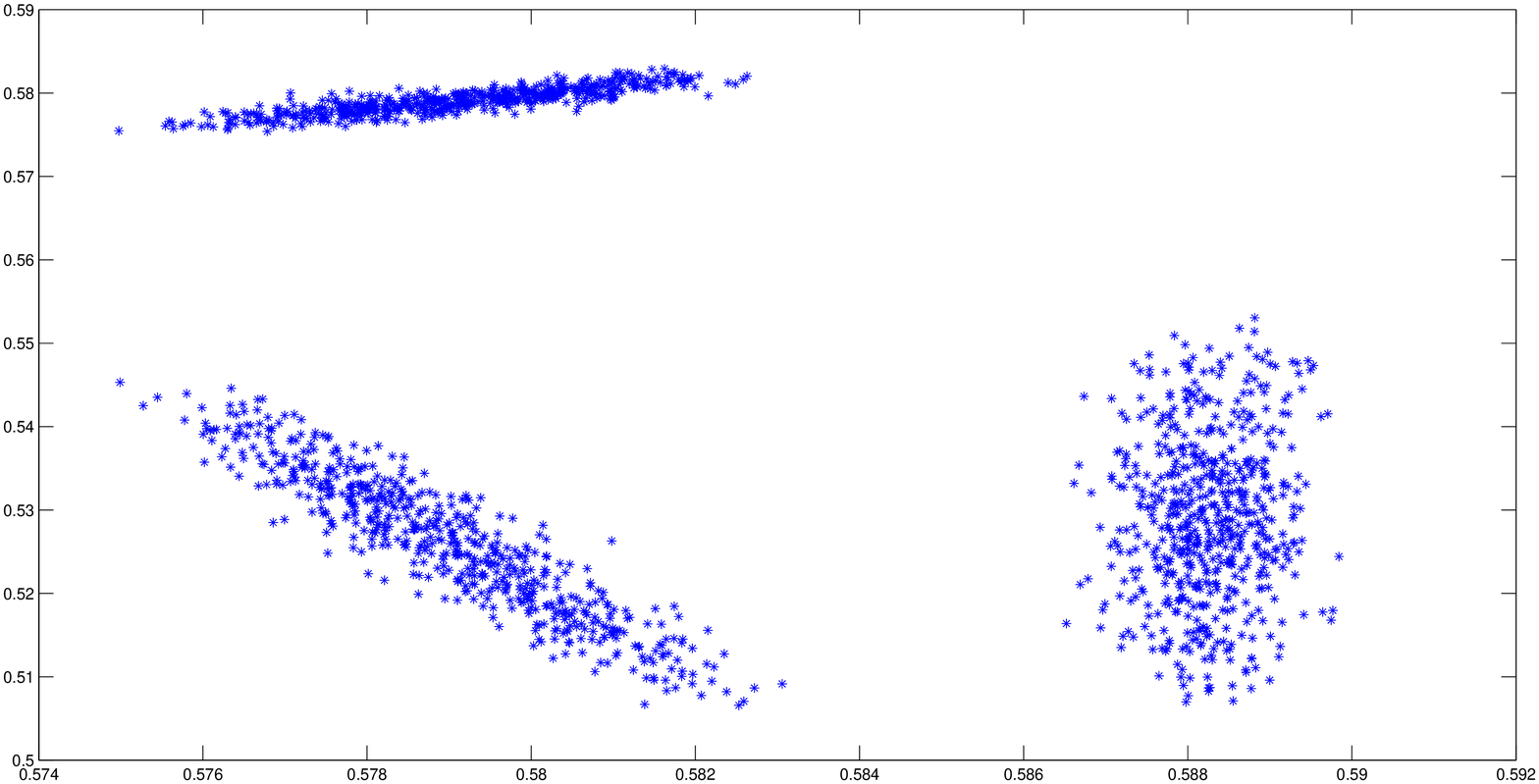}
\includegraphics[height=7cm,width=7.5cm]{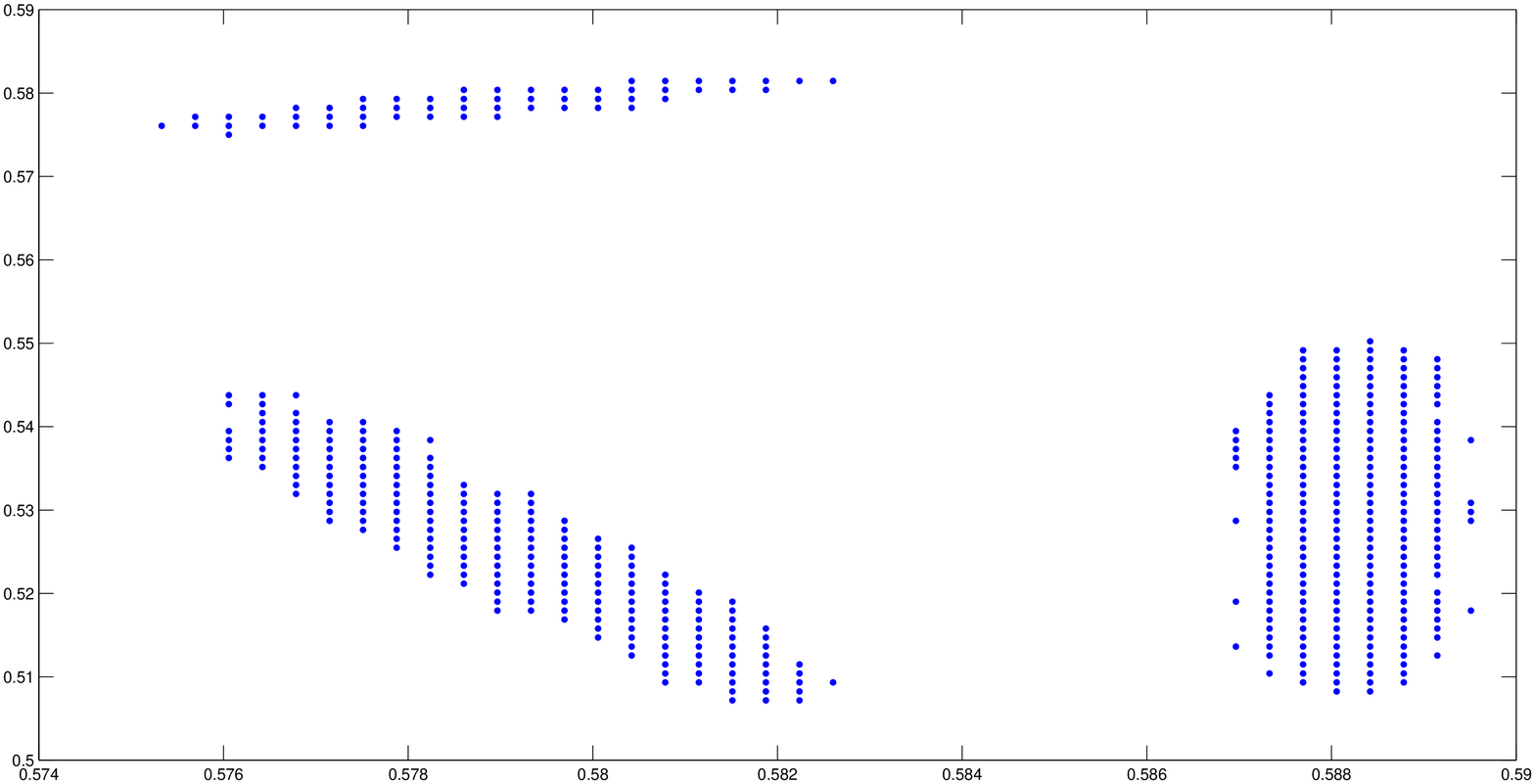}
\caption{The point cloud before (left) and after (right) denoising}
\label{denoise2}
\end{figure}

\begin{figure}
\includegraphics[height=7cm,width=7.5cm]{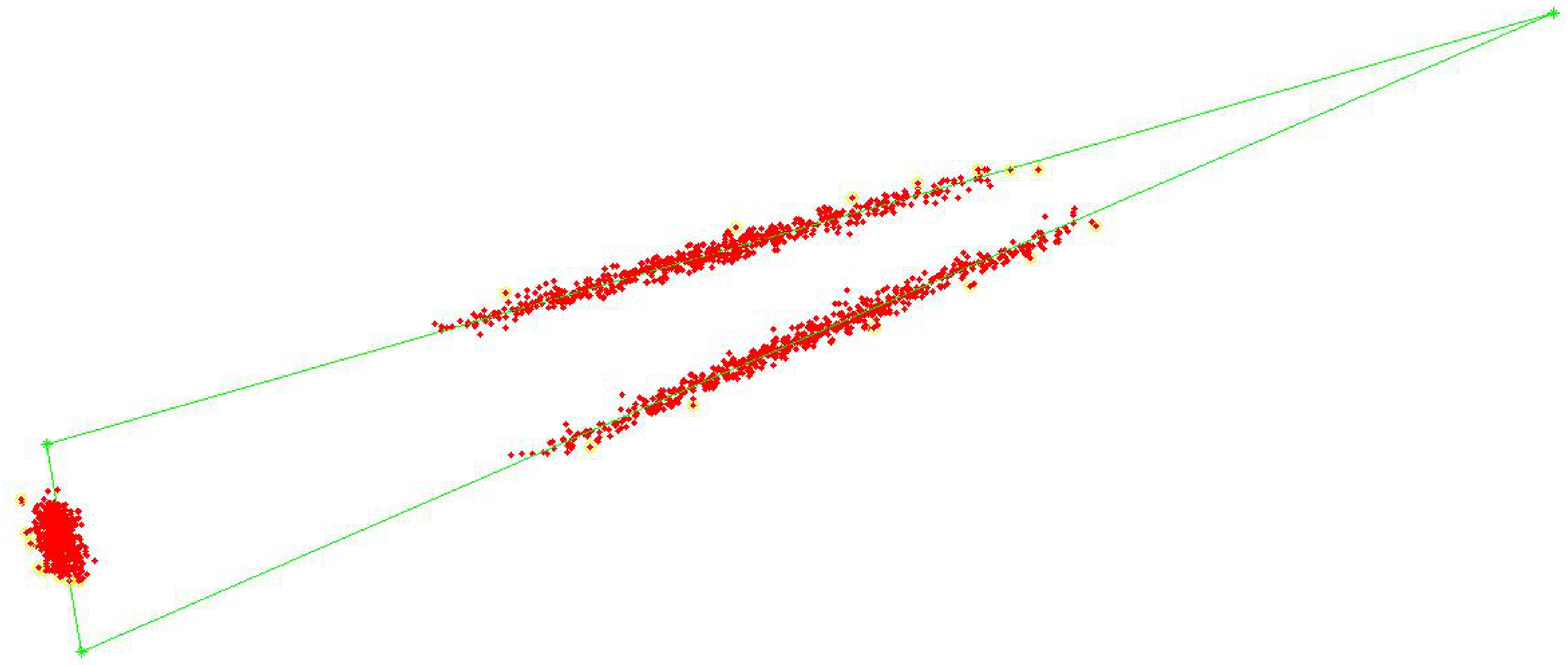}
\includegraphics[height=7cm,width=7.5cm]{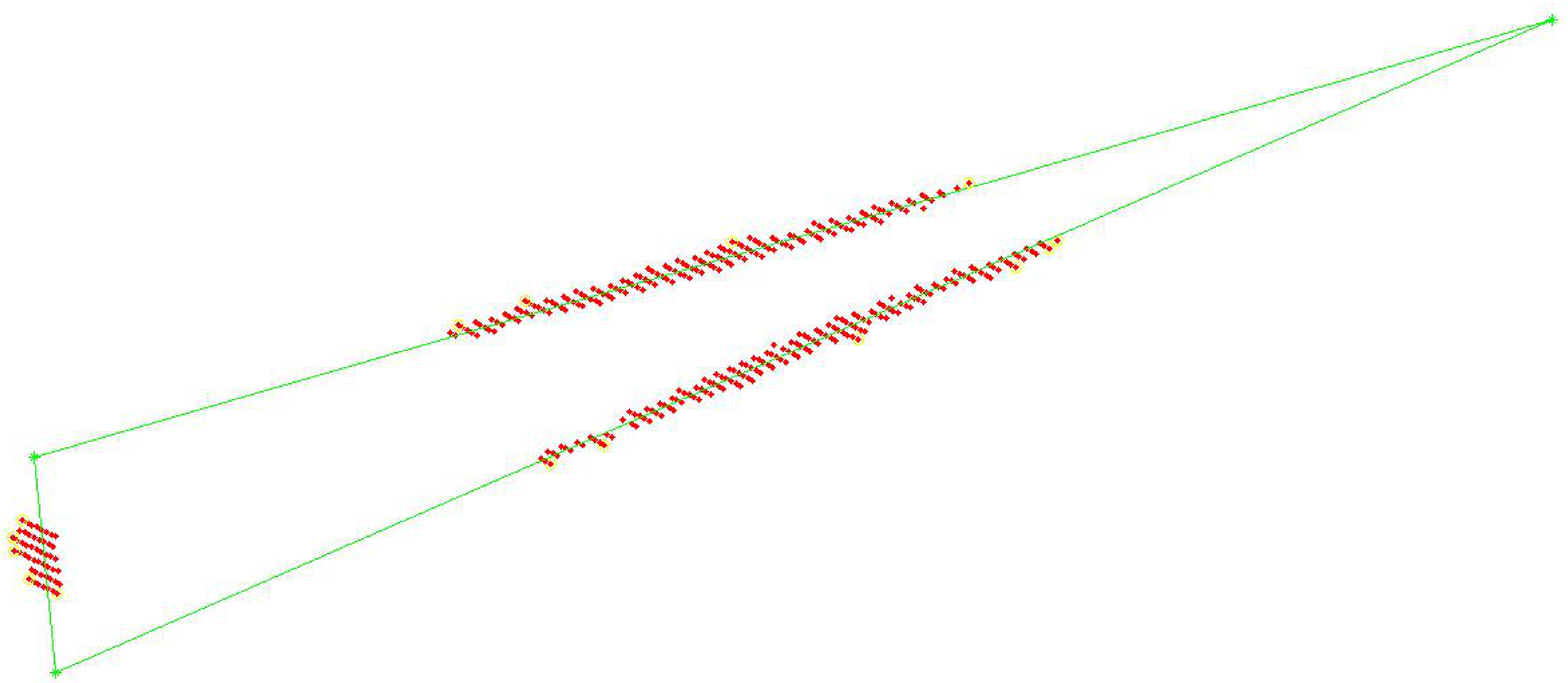}
\caption{Computational results of planes from noisy data (left) and denoised data (right).  The Green lines are the detected planes using the method proposed in the paper.}
\label{noiseVSdenoise}
\end{figure}

\begin{figure}
\includegraphics[height=6cm,width=15cm]{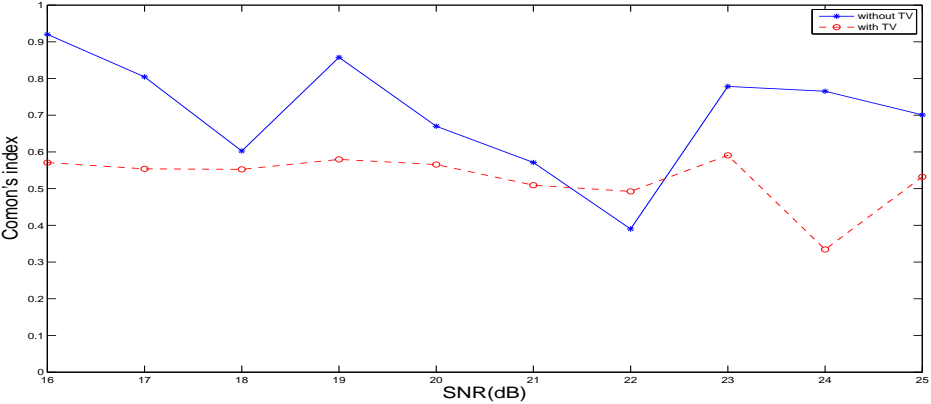}
\caption{Compare Comon's indices of FCA with and without TV denosing at high noise levels.}
\label{idx3}
\end{figure}

We proceed to solve the Rudin-Osher-Fatemi (ROF) model to obtain a denoised distance function $u(x,y)$
\begin{equation*}
\min_{u}\mathrm{TV}(u) + \lambda/2 || d - u ||^2_2,
\end{equation*}
where $TV(u)$ is the total variation of $u$ defined as
\begin{equation*}
TV(u) = \sum_{i,j}\sqrt{|u_{i+1,j} - u_{i,j}|^2 + |u_{i,j+1} - u_{i,j}|^2 }
\end{equation*}
for an image.  A variation for ease of minimization is:
\begin{equation*}
TV(u) = \sum_{i,j}|u_{i+i,j} - u_{i,j}| + |u_{i,j+1} - u_{i,j}|\;.
\end{equation*}
We shall use the recent {\em Chambolle's Algorithm} \cite{Chambolle_04} to solve this minimization problem.
Note that the smaller the parameter $\lambda$, the stronger the denoising.  Then the zero level set of the resulting minimizer $u(x,y)$ will be taken as the denoised point cloud.  In the real calculation, we will consider the following set with threshold
$\mathcal{{S}} = \{(x,y): u(x,y) = \tau\}$ where $\tau$ takes on a tiny value.  The noisy point cloud and the result after the noise removal are depicted in Fig. \ref{denoise2}.   The detected planes from both of them are shown in Fig. \ref{noiseVSdenoise}, and their intersections, i.e., the vertexes of the cone.  It can be noted that total variation denoising is very effective at preserving edges (thick lines in the figures) whilst smoothing away noise in flat regions.  The idea of denoising distance function by total variation extends to point cloud of any dimension.
We conduct experiments on performances of FCA with and without TV denoising were conducted.  A comparison of their performance is showed in Fig. \ref{idx3}. The mixtures are corrupted by additive white Gaussian noises varying from 16 dB to 25 dB. In general, TV denoising lowers Comon's indices at high noise level resulting in better separation.

\section{Concluding Remarks}
In this paper, we developed a novel BSS method based on facet component analysis.  We presented a facet based solvability condition
for the unique solvability of the nonnegative blind source separation problem up to scaling and permutation. Our approach exploited both the geometry of data matrix and the sparsity of the source signals.  Numerical results on NMR signals validated the solvability condition, and showed satisfactory performance of the proposed algorithm.
For noisy data, total variation denoising method serves as
a viable preprocessing step.

A line of future work is to separate more source signals from their mixtures,
known as an undetermined blind source separation, or uBSS.  This problem presents more challenge than the determined or over-determined BSS in that the mixing matrix is non-invertible. Some recent study has been done by two of the authors based on a geometric approach to retrieve the mixing matrix under suitable solvability conditions \cite{SX2}.

\end{document}